\documentclass[11pt,letter]{amsproc}
\usepackage{amssymb,amsmath,amsthm,anysize}
\usepackage{longtable}
\usepackage{booktabs}
\usepackage{multirow,bigstrut,bigdelim,color}
\usepackage{lscape}

\newtheorem{Thm}{Theorem}[section]

\newtheorem{Def}[Thm]{Definition}

\newtheorem{Pro}[Thm]{Proposition}
\newtheorem{La}[Thm]{Lemma}

\newtheorem{Claim}[Thm]{Claim}

\newenvironment{Prf}{\noindent\textbf{Proof.}}{\hfill $\Box$ \medskip}

\newcommand{\cV}{\mathcal{V}}
\newcommand{\cC}{\mathcal{C}}
\newcommand{\cP}{\mathcal{P}}
\newcommand{\F}{\mathbb{F}}

\newcommand{\GL}{\mathrm{GL}}

\newcommand{\PGL}{\mathrm{PGL}}

\newcommand{\PG}{\mathrm{PG}}

\renewcommand\le{\leqslant}
\renewcommand\ge{\geqslant}

\renewcommand{\#}{\times}

\begin{document}

\title{The symmetric representation of lines in ${\mathrm{PG}}(\F^3\otimes \F^3)$}
\author{Michel Lavrauw, Tomasz Popiel}

\address{
Michel Lavrauw, Sabanc{\i} University, Istanbul, Turkey 
\newline Email: \texttt{mlavrauw@sabanciuniv.edu}
\newline \indent Tomasz Popiel, Queen Mary University of London 
\newline Email: \texttt{t.popiel@qmul.ac.uk}}




\maketitle

\begin{abstract}
Let $\F$ be a finite field, an algebraically closed field, or the field of real numbers. 
Consider the vector space $V=\F^3 \otimes \F^3$ of $3 \times 3$ matrices over $\F$, and let $G \leq \text{PGL}(V)$ be the setwise stabiliser of the corresponding Segre variety $S_{3,3}(\F)$ in the projective space $\PG(V)$. 
The $G$-orbits of lines in $\text{PG}(V)$ were determined by the first author and Sheekey as part of their classification of tensors in $\F^2 \otimes V$ in the article ``Canonical forms of $2 \times 3 \times 3$ tensors over the real field, algebraically closed fields, and finite fields'', {\em Linear Algebra Appl.} {\bf 476} (2015) 133--147. 
Here we solve the related problem of classifying those line orbits that may be represented by {\em symmetric} matrices, or equivalently, of classifying the line orbits in the $\F$-span of the Veronese variety $\mathcal{V}_3(\F) \subset S_{3,3}(\F)$ under the natural action of $K=\PGL(3,\F)$. 
Interestingly, several of the $G$-orbits that have symmetric representatives split under the action of $K$, and in many cases this splitting depends on the characteristic of $\F$. Although our main focus is on the case where $\F$ is a finite field, our methods (which are mostly geometric) are easily adapted to include the case where $\F$ is an algebraically closed field, or the field of real numbers.
The corresponding orbit sizes and stabiliser subgroups of $K$ are also determined in the case where $\F$ is a finite field, and connections are drawn with old work of Jordan and Dickson on the classification of pencils of conics in $\PG(2,\F)$, or equivalently, of pairs of ternary quadratic forms over $\F$.
\end{abstract}

\section{Introduction} \label{intro}

\subsection{Set-up and summary of our results}
Consider the vector space $V=\F^3\otimes \F^3$ of $3\times 3$ matrices over a field $\F$, and recall that the corresponding {\em Segre variety} $S_{3,3}(\F)$ in the projective space $\PG(V) \cong \PG(8,\F)$ is the image of the map taking $(\langle v \rangle,\langle w \rangle) \in \PG(\F^3) \times \PG(\F^3)$ to $\langle v \otimes w \rangle$.
Let $G$ denote the setwise stabiliser of $S_{3,3}(\F)$ inside the projective general linear group $\PGL(V)$. 
The classification of $G$-orbits of lines in $\PG(V)$ was obtained by the first author and Sheekey~\cite{LaSh2015} as a consequence of their classification of tensors in $\F^2 \otimes \F^2\otimes \F^3$ for $\F$ a finite field, $\F$ an algebraically closed field, and $\F=\mathbb R$. 
This led to the classification \cite{LaSh2017} of all subspaces of $\PG(\F^2\otimes \F^3)$, and of the tensor orbits in $\F^2\otimes\F^3\otimes \F^r$ for every $r\geq 1$.

Here we study the {\em symmetric representation} of the line orbits in $\PG(V)$, by which we mean the following. 
Let $\mathcal{O}$ be a $G$-orbit of lines in ${\mathrm{PG}}(V)$, and consider the subspace $V_\text{s} \le V$ of symmetric $3 \times 3$ matrices over $\F$. 
If $\mathcal{O}$ happens to contain a line $L$ in $\PG(V_\text{s})$, then $L$ is called a {\it symmetric representative} of $\mathcal{O}$. 
If $\mathcal{O}$ has two symmetric representatives that are not in the same orbit under the natural action of $K=\PGL(3,\F)$, whereby a symmetric matrix $M$ is mapped by $D \in \GL(3,\F)$ to $DMD^\top$, then we say that the $G$-orbit $\mathcal{O}$ {\it splits} (under this action of $K$).

We address the following natural problems concerning the $G$-orbits of lines in $\PG(V)$:
\begin{itemize}
\item[(i)] We determine which $G$-orbits of lines in $\PG(V)$ have a symmetric representative.
\item[(ii)] We classify those orbits that have symmetric representatives, under the action of $K$.
\item[(iii)] In the case where $\F$ is a finite field, we determine for each $K$-orbit the corresponding stabiliser subgroup of $K$ and the orbit size.
\end{itemize}
Note that problem~(ii) is equivalent to the classification of $K$-orbits of lines in the $\F$-span $\langle \cV_3(\F) \rangle$ of the {\em Veronese variety}, or {\em quadric Veronesean}, $\cV_3(\F) \subset S_{3,3}(\F)$, namely the image of the {\em Veronese map} $\nu_3 : \PG(2,\F) \rightarrow \PG(5,\F)$ induced by the mapping taking $u \in \F^3$ to $u \otimes u$. 

Our main results are the solutions to problems~(i) and~(ii) for the case of a finite field $\F$. 
These are addressed in Section~\ref{sec:orbits} and summarised in Table \ref{table:main}. 
There are 14 orbits of lines in ${\mathrm{PG}}(V)$ under $G$, arising from the tensor orbits $o_4,\ldots,o_{17}$ in $\F^2 \otimes V$, in the notation of \cite{LaSh2015}, which we adopt here for consistency. 
Of these 14 orbits, only three do not have symmetric representatives, namely those arising from the tensor orbits $o_4$, $o_7$ and $o_{11}$. 
Moreover, the line orbits (corresponding to the tensor orbits) $o_5$, $o_6$, $o_9$, $o_{10}$ and $o_{17}$ do not split for any value of the characteristic $\operatorname{char}(\F)$ of $\F$, while $o_{14}$ and $o_{15}$ split for odd characteristic but not for even characteristic, $o_{12}$ and $o_{16}$ split for even characteristic but not for odd characteristic, and $o_8$ and $o_{13}$ split for all values of $\operatorname{char}(\F)$. 
We note that no $G$-line orbit splits into more than two $K$-line orbits. 
For algebraically closed fields and $\F=\mathbb R$, problems~(i) and~(ii) are handled in Section~\ref{sec:otherF}. 
The situation is overall somewhat simpler than in the finite case, but we note in particular that the results for an algebraically closed field $\F$ do depend on whether or not $\operatorname{char}(\F)=2$.

Problem~(iii) is addressed in Section~\ref{sec:stabs}, with the results summarised in Table~\ref{table:stabs}. 
As noted above, the splitting (or not) of a $G$-orbit under $K$ sometimes depends on whether $\operatorname{char}(\F)=2$ or not. 
Moreover, the structures of the corresponding line stabilisers inside $K$ can also depend on $\operatorname{char}(\F)$. 
It seems remarkable, therefore, that the number of symmetric representatives of any given $G$-orbit turns out to be independent of $\operatorname{char}(\F)$ (see Table~\ref{table:orbitsLengths}). 

\subsection{Historical context and commentary on our results} \label{sec:history}
Because our results imply the classification of lines in $\langle \cV_3(\F) \rangle$ under the natural action of $\PGL(3,\F)$, when $\operatorname{char}(\F)\neq 2$ they also imply the classification of {\em pencils of conics} in $\PG(2,\F)$, namely, one-dimensional subspaces of ternary quadratic forms over $\F$. 
We refer the reader to Section~\ref{sec:pencils} for details about this correspondence.
The latter classification problem goes back to old work of Jordan~\cite{Jordan1906,Jordan1907} and Dickson~\cite{Dickson1908}. 
The classification over the reals and the complex numbers was obtained by Jordan  \cite{Jordan1906,Jordan1907} in 1906--1907; there are 13 and 8 orbits, respectively (in accordance with our results in Section~\ref{sec:otherF}). 
This classification was later extended to algebraically closed fields using the theory of matrix elementary divisors (due to Weierstrass). 
For example, Wall~\cite{Wall1977} refers to Segre's classification of pencils of quadrics, pointing out that it was well known at the time Wall's paper was published (1977), and remarking that it appears in various standard textbooks, including those of Gantmacher~\cite{Gantmacher} and Hodge and Pedoe~\cite{HodgePedoe}. 
However, as we explain below, this approach does not adequately treat the finite field case. 
Another approach which seems to work exclusively over the complex numbers is that due to Artamkin and Nurmiev \cite{ArNu2002}, who appeal to connections with the theory of Lie algebras.

The history of the problem of classifying pencils of conics over {\em finite} fields seems to be somewhat more complicated. 
In the odd characteristic case, the classification was obtained by Dickson~\cite{Dickson1908} in 1908. 
The even characteristic case was studied two decades later by Campbell~\cite{Campbell1927}, who provided a list of inequivalent classes of pencils of conics in $\PG(2,\F_q)$, $q$ even. 
However, unlike Dickson, Campbell did not obtain a full classification. 
Campbell was aware of this, stating on the first page of his paper that: {\it ``If there is an arbitrary coefficient in the typical pencil we say this pencil represents a set of classes, whenever different values of this coefficient may give nonequivalent pencils and so represent distinct classes.''}.
These ``sets of classes'' are listed on \cite[p.~406]{Campbell1927} as Set~10, Set~14, Set~15, Set~16 and Set~17.
This also explains why Campbell's paper is so short: the main difficulties in what would have been be a complete classification are not addressed. In particular, the pencils without binary forms, which in the odd characteristic case correspond to our aforementioned case $o_{17}$, are not classified (see also Dickson's comment from his paper \cite{Dickson1908}, quoted below). 
Although a full classification in the even characteristic case is alluded to in the literature, we have not been able to find an explicit list of orbits with an accompanying proof anywhere. 
In particular, Hirschfeld states the classification in full as Theorem 7.31 of his book \cite{Hirschfeld1998}, but attributes the result to Campbell \cite{Campbell1927}, who, as explained above, neither stated nor proved the complete classification. 
We intend to complete the classification in a forthcoming paper.

We now explain why the elementary divisor method used for algebraically closed fields by C.~Segre and others (as explained above) is inadequate to treat the finite field case. 
Hodge and Pedoe work over an algebraically closed field $\mathbb{F}$ and prove \cite[Theorem~I of Chapter~XIII]{HodgePedoe} that two $n \times n$ matrices $A$ and $B$, with $B$ non-singular, can be simultaneously transformed via a change of coordinates to matrices $C$ and $D$, with $D$ non-singular, respectively, if and only if the linear combinations $A-\lambda B$ and $C-\lambda D$ (where $\lambda$ is a variable) have the same elementary divisors. 
Note here that, when $\operatorname{char}(\F)$ is odd, a change of coordinates of the quadratic form associated with the matrix $A$, say, corresponds precisely to a mapping of $A$ to $ZAZ^\top$ for some non-singular matrix $Z$, and hence to the natural action of $\text{PGL}(n,\mathbb{F})$. 
This result is certainly false when ${\mathbb{F}}$ is finite.
For example, if we take $\mathbb{F}=\mathbb{F}_q$, $q$ even, and
\[
B=D=\left[ \begin{matrix} 0 & 0 & 1 \\ 0 & 1 & 0 \\ 1 & 0 & 0 \end{matrix} \right], \quad 
A=\left[ \begin{matrix} 0 & 0 & 0 \\ 0 & 0 & 1 \\ 0 & 1 & 0 \end{matrix} \right], \quad 
C=\left[ \begin{matrix} 0 & 0 & 0 \\ 0 & 0 & 1 \\ 0 & 1 & 1 \end{matrix} \right],
\]
then the matrices $A-\lambda B$ and $C-\lambda D$ have the same elementary divisors, namely $1$, $1$ and $\lambda^3$. 
However, there is no non-singular matrix $Z=(z_{ij})$ such that $ZAZ^\top=C$, because the $(3,3)$-entry of $ZAZ^\top$ is $2z_{32}z_{33}$, which is $0$ and hence not equal to the $(3,3)$-entry of $C$ in characteristic~$2$. 
This example, which in fact illustrates that \cite[Theorem~I of Chapter~XIII]{HodgePedoe} is false when $\operatorname{char}(\mathbb{F})=2$, arises from the $o_{16}$ case of our classification: by Table~\ref{table:main}, the line orbits represented by the matrices $Ax+By$ and $Cx+Dy$, where $x$ and $y$ are variables, are inequivalent under the natural action of $\text{PGL}(3,\mathbb{F}_q)$ when $q$ is even.
It is also straightforward to find a counterexample for finite fields of odd characteristic; in particular, one arises from the $o_8$ case of our classification. 
Hence, \cite[Theorem~I of Chapter~XIII]{HodgePedoe} does not imply a classification of lines in $\PG(\F^3\otimes \F^3)$ over finite fields.
Gantmacher~\cite{Gantmacher} takes a more refined approach, using both finite and infinite elementary divisors (due to Kronecker), but this approach is still inadequate for finite fields (for similar reasons).


We also remark on (what we feel are) some advantages of our approach as compared with Dickson's original proof, in the case of a finite field of {\em odd} characteristic. 
Dickson~\cite{Dickson1908} determined an exhaustive list of 15 equivalence classes of pairs of ternary quadratic forms over a finite field $\mathbb{F}_q$ with $q$ odd, consistent with our results in Table~\ref{table:main}. 
In the first paragraph of his paper, he anticipated that: {\it ``The main difficulty lies in the case in which the family contains no binary forms, and that in which the binary forms are all irreducible. Neither of these cases occur when the field is ${\mathbb{C}}$ or ${\mathbb{R}}$, so that the problem is quite simple for these fields.''} 
Indeed, 15 of the 18 pages of Dickson's paper \cite{Dickson1908} are dedicated to the classification of these two cases. 
The proof gives explicit coordinate transformations in order to reduce the families of quadratic forms to canonical representatives of the associated equivalence classes, and can at times be quite tedious. 
In particular, the cases $q \equiv 0$, $1$ and $2 \pmod 3$ are treated separately in the proof of the case in which the family contains no binary forms. 
This case also relies on knowledge of the number of irreducible cubics of a given form and refers to Dickson's treatise on linear groups \cite{Dickson????}. 
Our proof is quite different, and in particular we do not need to treat the cases $q \equiv 0$, $1$ and $2 \pmod 3$ separately. 
In fact, our approach applies more or less uniformly for both even and odd characteristic. 
(Moreover, as noted above, we compute stabilisers for all orbits, which do not seem to have been recorded anywhere as far as we can tell.) 
We note, however, that despite the differences between our approach and Dickson's, there are cases in which we encounter similar difficulties. 
In particular, in the case of pencils without binary forms, which correspond to lines in $\PG(\F_q^3\otimes \F_q^3)$ without points of rank~$3$ (namely those of type $o_{17}$, treated at the end of Section~\ref{sec:orbits}), both Dickson's proof and our proof are based on counting arguments. 
This seems unavoidable. 
Interestingly enough, the proof in \cite[Section 3.3]{LaSh2015} of the fact that there is a single orbit of constant rank-3 lines in $\PG(\F_q^3\otimes \F_q^3)$ is {\em also}, seemingly unavoidably, based on a counting argument. 
All three arguments are, however, counting different objects. 

Finally, we remark that our proofs are largely geometric in nature, in contrast with those of Dickson~\cite{Dickson1908}, and indeed with the aforementioned matrix elementary divisor arguments. 
This geometric approach has been our starting point towards a classification of orbits of {\em planes} in ${\mathrm{PG}}(\F^3\otimes \F^3)$ for $\F$ a finite field, which correspond to {\em nets} of conics, namely, 
{\em three}-dimensional subspaces of ternary quadratic forms.
This natural (and more complicated) next case was investigated by Wilson \cite{Wilson1914} and Campbell \cite{Campbell1928}, but a complete classification is still unknown. 
(We also remark that the case $\mathbb{F}=\mathbb{R}$ was treated by Wall~\cite{Wall1977}.)


\section{Preliminaries} \label{prelims}

Here we collect some preliminary information for background and later reference.

\subsection{Orbits of tensors in $\mathbb{F}^2 \otimes \mathbb{F}^3 \otimes \mathbb{F}^3$}

Write 
\[
V_1 = \F^2, \quad V=\F^3 \otimes \F^3 \quad \text{and} \quad \overline{V} = V_1 \otimes V, 
\]
and let $G$ be the setwise stabiliser in $\text{GL}(\overline{V})$ of the set of fundamental tensors in $\overline{V}$, namely the tensors of the form $v_1 \otimes v_2 \otimes v_3$ with $v_1 \in \F^2$ and $v_2,v_3 \in \F^3$. 
The $G$-orbits of tensors in $\overline{V}$ were classified in \cite[Main Theorem]{LaSh2015}.
In particular, in the case where $\F$ is a finite field $\F_q$, there are precisely $18$ orbits, with representatives given in terms of a basis $\{e_1,e_2,e_3\}$ for $\mathbb{F}_q^3$ in the table on \cite[p.~146]{LaSh2015}. 
For convenience, the information in that table is included here in Table~\ref{pointTable}.

\begin{table}[!t]
\begin{tabular}{llcl}
\toprule
Orbit & Representative & Condition & Rank dist. \\
\midrule
$o_0$ & $0$ && $[0,0,0]$ \\
$o_1$ & $e_1 \otimes e_1 \otimes e_1 $ && $[1,0,0]$ \\
$o_2$ & $e_1 \otimes (e_1 \otimes e_1+ e_2\otimes e_2)$ && $[0,1,0]$ \\
$o_3$ & $e_1 \otimes e$ && $[0,0,1]$ \\
$o_4$ & $e_1 \otimes e_1 \otimes e_1 + e_2\otimes e_1 \otimes e_2$ && $[q+1,0,0]$ \\
$o_5$ & $e_1 \otimes e_1 \otimes e_1 + e_2\otimes e_2 \otimes e_2$ && $[2,q-1,0]$ \\
$o_6$ & $e_1 \otimes e_1 \otimes e_1 + e_2\otimes (e_1 \otimes e_2 + e_2 \otimes e_1)$ && $[1,q,0]$ \\
$o_7$ & $e_1 \otimes e_1 \otimes e_3 + e_2\otimes (e_1 \otimes e_1 + e_2 \otimes e_2)$ && $[1,q,0]$ \\
$o_8$ & $e_1 \otimes e_1 \otimes e_1 + e_2\otimes (e_2 \otimes e_2 + e_3 \otimes e_3)$ && $[1,1,q-1]$ \\
$o_9$ & $e_1 \otimes e_3 \otimes e_1 + e_2\otimes e$ && $[1,0,q]$ \\
$o_{10}$ & $e_1\otimes (e_1\otimes e_1+ e_2\otimes e_2+u e_1\otimes e_2)+e_2\otimes (e_1\otimes e_2+v e_2\otimes e_1)$ & ($*$) & $[0,q+1,0]$ \\
$o_{11}$ & $e_1\otimes (e_1 \otimes e_1 + e_2 \otimes e_2) + e_2\otimes (e_1 \otimes e_2 + e_2 \otimes e_3)$ && $[0,q+1,0]$ \\
$o_{12}$ & $e_1\otimes (e_1 \otimes e_1 + e_2 \otimes e_2) + e_2\otimes (e_1 \otimes e_3 + e_3 \otimes e_2)$ && $[0,q+1,0]$ \\
$o_{13}$ & $e_1\otimes (e_1 \otimes e_1 + e_2 \otimes e_2) + e_2\otimes (e_1 \otimes e_2 + e_3 \otimes e_3)$ &&$[0,2,q-1]$ \\
$o_{14}$ & $e_1\otimes (e_1 \otimes e_1 + e_2 \otimes e_2) + e_2\otimes (e_2 \otimes e_2 + e_3 \otimes e_3)$ && $[0,3,q-2]$ \\
$o_{15}$ & $e_1\otimes (e+u e_1\otimes e_2) + e_2\otimes (e_1\otimes e_2+v e_2\otimes e_1)$ & ($*$) & $[0,1,q]$ \\
$o_{16}$ & $e_1\otimes e + e_2\otimes (e_1 \otimes e_2 + e_2 \otimes e_3)$ && $[0,1,q]$ \\
$o_{17}$ & $e_1\otimes e + e_2\otimes (e_1\otimes e_2 + e_2\otimes e_3 + e_3\otimes (\alpha e_1 + \beta e_2 + \gamma e_3))$ & ($**$) & $[0,0,q+1]$ \\
\bottomrule
\end{tabular}
\caption{Orbits of tensors in $\overline{V} = \mathbb{F}_q^2 \otimes \mathbb{F}_q^3 \otimes \mathbb{F}_q^3$ under the setwise stabiliser in $\text{GL}(\overline{V})$ of the set of fundamental tensors in $\overline{V}$, as per \cite[p.~146]{LaSh2015}. 
Representatives are given in terms of a basis $\{e_1,e_2,e_3\}$ of $\mathbb{F}_q^3$, with $e = \sum_{i=1}^3 e_i \otimes e_i$. 
The final column shows the rank distribution of the first contraction space of each representative.
Condition~($*$) is: $v\lambda^2+uv\lambda - 1 \neq 0$ for all $\lambda \in \F_q$. 
Condition~($**$) is: $\lambda^3+\gamma \lambda^2- \beta \lambda+ \alpha \neq 0$ for all $\lambda \in \F_q$.}
\label{pointTable}
\end{table}

In this paper, we are interested in symmetric representatives of line orbits in the projective space $\PG(V)$. 
The line orbits themselves can be obtained by considering the {\em first contraction spaces} of tensors in $\overline{V}$. 
As per \cite[p.~136]{LaSh2015}, the {\em first contraction space} of a tensor $A \in \overline{V}$ is the subspace
\[
A_1 = \langle w_1^\vee(A) : w_1^\vee \in V_1^\vee \rangle
\]
of $V$. 
Here $V_1^\vee$ is the dual of $V_1$, and $w_1^\vee(A)$ is defined by its action on fundamental tensors via $w_1^\vee(v_1 \otimes v_2 \otimes v_3) = w_1^\vee(v_1)v_2 \otimes v_3$. 
Recall also that the {\em rank} of a point in $\PG(V)$ is the rank of any ($3 \times 3$) matrix representing that point (and that this does not depend on the choice of bases for the factors of the tensor product). 
In geometric terms, a point has rank~1 if it is contained in the Segre variety $S = S_{3,3}(\F) \subset \PG(V)$, rank~2 if is not contained in $S$ but is contained in the secant variety of $S$, and rank~3 if it is not contained in the secant variety of $S$. 
Note that the last column of Table~\ref{pointTable} shows, for each orbit, the {\em rank distribution} of the first contraction space of a representative $A$, namely a list $[a_1,a_2,a_3]$ where $a_i$ is the number of points of rank $i$ in $\text{PG}(A_1)$. 
The first contraction spaces of the tensors in orbits $o_4,\ldots,o_{17}$ are lines of $\PG(V)$; in particular, their rank distributions satisfy $a_1+a_2+a_3=q+1$.

\subsection{Properties of the quadric Veronesean} \label{V3props}

Here we collect some facts about the quadric Veronesean $\cV_3(\F)$ which are used throughout the paper. 
Most of these properties belong to the folklore of classical algebraic geometry. For proofs and/or further details, the reader may consult a standard reference such as Harris~\cite[p.~23]{Harris} or Hirschfeld and Thas~\cite[Chapter~4]{HiTh2016}. We restrict ourselves to those results that are frequently used in our proofs. 

The quadric Veronesean has been studied for over a century and many interesting results have been obtained.
A remarkable characterisation of $\cV_3({\mathbb{C}})$ was, for instance, given by Bertini in 1923 \cite{Bertini1923}. This was extended by Mazzocca and Melone to $\cV_3(\F_q)$. Their paper \cite{MaMe1984} also contains an interesting list of references to some of the earlier works on the subject by Italian geometers. For a longer list of references we refer the reader to the more general survey by Havlicek~\cite{Havlicek2003} on Veronese varieties over fields of positive characteristic.

As noted in Section~\ref{intro}, $\cV_3(\F)$ is the image of the map $\nu_3 : \PG(2,\F) \rightarrow \PG(5,\F)$ induced by the mapping taking $u \in \F^3$ to $u \otimes u$. 
The subgroup $K$ of $\PGL(6,\F)$, isomorphic to $\PGL(3,\F)$, with $D \in \GL(3,\F)$ mapping a symmetric matrix $M$ to $M^D = DMD^\top$, is equal to the setwise stabiliser
of $\cV_3(\F)$ unless $\F=\F_2$, in which case $K$ is a proper subgroup of the setwise stabiliser (see e.g. \cite[p.~148]{HiTh2016}).
We also record the following facts, which are readily obtained from the relevant definitions:

\begin{itemize}
\item[(F1)] The image of a line of $\PG(2,\F)$ under $\nu_3$ is a conic. 
A plane of $\PG(5,\F)$ intersecting $\cV_3(\F)$ in the image of line of $\PG(2,\F)$ is called a {\it conic plane}. 
Unless $\F=\F_2$, each plane intersecting $\cV_3(\F)$ in a conic is a conic plane, while for $\F=\F_2$ a conic consists only of three points, and there are planes intersecting $\cV_3(\F_2)$ in three points which are not conic planes. 
Note that this definition of conic planes of $\cV_3(\F)$ is consistent with \cite[p.~148]{HiTh2016}. Similarly, we define {\it conics in $\cV_3(\F)$} as images of lines of $\PG(2,\F)$.
\item[(F2)] Each two points $P,Q$ of $\cV_3(\F)$ lie on a unique conic $\cC(P,Q)$ in $\cV_3(\F)$, given by $\cC(P,Q)=\nu_3(\langle \nu_3^{-1}(P),\nu_3^{-1}(Q)\rangle)$.
\item[(F3)] Each rank-2 point $R$ in $\langle \cV_3(\F) \rangle$ determines a unique conic $\cC(R)$ in $\cV_3(\F)$.
The point $R$ is called an {\it exterior point} if it lies on a tangent to $\cC(R)$, and an {\it interior point} otherwise. When $\operatorname{char}(\F)=2$, there are no interior points, because all the tangent lines of a conic are concurrent; their common point is called the {\it nucleus} of the conic.
\item[(F4)] The quadrics of $\PG(2,\F)$ are mapped by $\nu_3$ onto the hyperplane sections of $\cV_3(\F)$. 
A conic consisting of just one point (two distinct lines over the quadratic extension) corresponds to a hyperplane intersecting $\cV_3(\F)$ in one point. 
A repeated line of $\PG(2,\F)$ corresponds to a hyperplane meeting $\cV_3(\F)$ in a conic; two distinct lines correspond to a hyperplane meeting $\cV_3(\F)$ in two conics; and a non-degenerate conic corresponds to a hyperplane meeting $\cV_3(\F)$ in a normal rational curve.
\item[(F5)] If $\operatorname{char}(\F)=2$ then the nuclei of all of the conics contained in $\cV_3(\F)$ form a plane, called the {\em nucleus plane} of $\cV_3(\F)$. 
In the representation of the points of $\cV_3(\F)$ as symmetric $3\times 3$ matrices of rank~1, the nucleus plane comprises the matrices with zeroes on the main diagonal (with no restriction on the other three variables).
\end{itemize}

The $K$-orbits of {\em points} in $\langle \cV_3(\F)\rangle$ are well understood. 
For convenience, we note some facts about these point orbits in the case where $\F$ is a finite field $\F_q$:
\begin{itemize}
\item There is one $K$-orbit of points of rank~1: $K$ acts transitively on the set of 4-tuples of points of $\cV_3(\F_q)$, no three of which are on a conic. 
This orbit has size $q^2+q+1$ (the number of points in $\PG(2,\F_q)$).
\item There are two $K$-orbits of points of rank~2. 
In odd characteristic, one of these orbits consists of all the exterior points, and the other consists of all the interior points.
Denoting these orbits by $\cP_{2,\text{e}}$ and $\cP_{2,\text{i}}$ respectively, we have 
$|\cP_{2,\text{e}}|=\frac{1}{2}q(q+1)(q^2+q+1)$ and $|\cP_{2,\text{i}}|=\frac{1}{2}q(q-1)(q^2+q+1)$. 
In even characteristic, one orbit consists of all the points that lie on the nucleus plane of $\cV_3(\F_q)$, and the other orbit consists of all the other points of rank~2. 
Denoting these orbits by $\cP_{2,\text{n}}$ and $\cP_{2,\text{s}}$ respectively, we have 
$|\cP_{2,\text{n}}|=q^2+q+1$ and $|\cP_{2,\text{s}}|=(q^2-1)(q^2+q+1)$. 
Hence, regardless of the value of $q$, the total number of points of rank~2 is $q^2(q^2+q+1)$.
\item Finally, there is a unique $K$-orbit of rank~3 points, of size
\[
\frac{q^6-1}{q-1}-(q^2+1)(q^2+q+1)=q^5-q^2.
\]
\end{itemize}

\section{Line orbits in $\langle \mathcal{V}_3(\F) \rangle$ for $\F$ a finite field}\label{sec:orbits}

We now address problems~(i) and~(ii) of Section~\ref{intro} in the case where $\F$ is a finite field $\F_q$. 
Our strategy is as follows. 
We consider each of the orbits $o_0,\ldots,o_{17}$ of tensors in $\F_q^2 \otimes \F_q^3 \otimes \F_q^3$, which are shown in Table~\ref{pointTable}.  
Given a representative of an orbit $o_i$ from the second column of Table~\ref{pointTable}, we consider the corresponding first contraction space $M_i$, which is a subspace of $V=\F_q^3 \otimes \F_q^3$. 
When $i \ge 4$ in Table~\ref{pointTable}, $\text{PG}(M_i)$ is a line of $\PG(V)$, comprising $q+1$ points, which we represent by $3 \times 3$ matrices. 
For each $i \ge 4$, we first address problem~(i) by checking whether $\text{PG}(M_i)$ can be mapped into $\langle \mathcal{V}_3(\F_q) \rangle$ by the action of $\text{PGL}(3,\F_q) \times \text{PGL}(3,\F_q)$ induced by the action of $\GL(3,\F_q) \times \GL(3,\F_q)$ taking a $3 \times 3$ matrix $M$ to $M^{(B,C)} = BMC$ (where $B,C \in \GL(3,\F_q)$). 
In other words, we check whether $M_i^{(B,C)}$ can be a subspace of {\em symmetric} $3 \times 3$ matrices. 
If it cannot, then the $G$-line orbit arising from the tensor orbit $o_i$ is not represented in $\langle \mathcal{V}_3(\F_q) \rangle$, that is, it does not have a symmetric representative in the sense defined in Section~\ref{intro}.
If it can, then we address problem~(ii) by determining the orbits of the group $K = \text{PGL}(3,\F_q)$ in the action $M^D = DMD^\top$ for $M$ a symmetric matrix and $D \in \GL(3,\F_q)$.

Note also that when considering problem~(i) as described above, we may take $C$ to be the identity matrix, because the line $\PG(BM_iC)$ is equivalent under the action of $K$ to the projective space obtained from the vector subspace
\[
(BM_iC)^{(C^{-1})^\top} = (C^{-1})^\top(BM_iC)C^{-1} = (C^{-1})^\top B M_i.
\] 
This simplifies the proof of the fact that certain tensor orbits, namely $o_4$, $o_7$ and $o_{11}$, do {\em not} yield lines with symmetric representatives. 
The remaining tensor orbits $o_i$ (with $i \ge 4$) do yield lines with symmetric representatives, and the representatives of the corresponding $K$-orbits are listed in Table~\ref{table:main}.
It turns out that in each case there are at most two $K$-orbits of lines. 
The $K$-orbit in the second column of the table arises for all values of $q$, and sometimes there is another $K$-orbit, with representative shown in the third column if $q$ is odd and in the fourth column if $q$ is even. 
The following notation is used for brevity in Table~\ref{table:main} (and in the proofs):

\begin{Def} \label{def:notation}
\textnormal{
The matrices in Table~\ref{table:main} represent subspaces of symmetric matrices over the finite field $\mathbb{F}_q$. 
The subscript ``$x,y$'' indicates that the pair $(x,y)$ ranges over all values in $\mathbb{F}_q^2$, and the symbol $\cdot$ denotes $0$. 
For example, in the first line of the table,
\[
\left[ \begin{matrix} x & \cdot & \cdot \\ \cdot & y & \cdot \\ \cdot & \cdot & \cdot \end{matrix} \right]_{x,y} = 
\left\{ \left[ \begin{matrix} x & 0 & 0 \\ 0 & y & 0 \\ 0 & 0 & 0 \end{matrix} \right] : (x,y) \in \mathbb{F}_q^2 \right\},
\]
and the line orbit representative in $\langle \cV_3(\F_q) \rangle$ is the corresponding projective space.
The symbol $\Box$ is used to denote the set of squares in $\F_q$.
}
\end{Def}

\begin{table}[!t]
\begin{tabular}{ccccc}
\toprule
Tensor & \multicolumn{3}{c}{Line orbit representatives in $\langle \mathcal{V}_3(\mathbb{F}_q) \rangle$} & Conditions \\
\cline{2-4} orbit & Common orbit (all $q$) & \multicolumn{2}{c}{Additional orbit} & \\
\cline{3-4} && $q$ odd & $q$ even & \\
\midrule
$o_5$ & $\left[ \begin{matrix} x & \cdot & \cdot \\ \cdot & y & \cdot \\ \cdot & \cdot & \cdot \end{matrix} \right]_{x,y}$ &&& \\
\addlinespace[2pt]
$o_6$ & $\left[ \begin{matrix} x & y & \cdot \\ y & \cdot & \cdot \\ \cdot & \cdot & \cdot \end{matrix} \right]_{x,y}$ &&& \\
\addlinespace[2pt]
$o_8$ & $\left[ \begin{matrix} x & \cdot & \cdot \\ \cdot & y & \cdot \\ \cdot & \cdot & y \end{matrix} \right]_{x,y}$ 
& $\left[ \begin{matrix} x & \cdot & \cdot \\ \cdot & y & \cdot \\ \cdot & \cdot & \gamma y \end{matrix} \right]_{x,y}$
& $\left[ \begin{matrix} x & \cdot & \cdot \\ \cdot & \cdot & y \\ \cdot & y & \cdot \end{matrix} \right]_{x,y}$ 
& $\gamma \not \in \Box$ \\
\addlinespace[2pt]
$o_9$ & $\left[ \begin{matrix} x & \cdot & y \\ \cdot & y & \cdot \\ y & \cdot & \cdot \end{matrix} \right]_{x,y}$ &&& \\
\addlinespace[2pt]
$o_{10}$ & $\left[ \begin{matrix} vx & y & \cdot \\ y & x+uy & \cdot \\ \cdot & \cdot & \cdot \end{matrix} \right]_{x,y}$ &&& ($*$) \\
\addlinespace[2pt]
$o_{12}$ & $\left[ \begin{matrix} \cdot & x & \cdot \\ x & \cdot & y \\ \cdot & y & \cdot \end{matrix} \right]_{x,y}$ & 
& $\left[ \begin{matrix} \cdot & x & \cdot \\ x & x+y & y \\ \cdot & y & \cdot \end{matrix} \right]_{x,y}$ & \\
\addlinespace[2pt]
$o_{13}$ & $\left[ \begin{matrix} \cdot & x & \cdot \\ x & y & \cdot \\ \cdot & \cdot & y \end{matrix} \right]_{x,y}$ 
& $\left[ \begin{matrix} \cdot & x & \cdot \\ x & y & \cdot \\ \cdot & \cdot & \gamma y \end{matrix} \right]_{x,y}$
& $\left[ \begin{matrix} \cdot & x & \cdot \\ x & x+y & \cdot \\ \cdot & \cdot & y \end{matrix} \right]_{x,y}$ 
& $\gamma \not \in \Box$ \\
\addlinespace[2pt]
$o_{14}$ & $\left[ \begin{matrix} x & \cdot & \cdot \\ \cdot & x+y & \cdot \\ \cdot & \cdot & y \end{matrix} \right]_{x,y}$ 
& $\left[ \begin{matrix} x & \cdot & \cdot \\ \cdot & \gamma(x+y) & \cdot \\ \cdot & \cdot & y \end{matrix} \right]_{x,y}$
&& $\gamma \not \in \Box$ \\
\addlinespace[2pt]
$o_{15}$ & $\left[ \begin{matrix} v_1y & x & \cdot \\ x & ux+y & \cdot \\ \cdot & \cdot & x \end{matrix} \right]_{x,y}$ 
& $\left[ \begin{matrix} v_2y & x & \cdot \\ x & ux+y & \cdot \\ \cdot & \cdot & x \end{matrix} \right]_{x,y}$
&& ($*$), $\begin{array}{ll}-v_1 \in \Box \\ -v_2 \not \in \Box \end{array}$ \\
\addlinespace[2pt]
$o_{16}$ & $\left[ \begin{matrix} \cdot & \cdot & x \\ \cdot & x & y \\ x & y & \cdot \end{matrix} \right]_{x,y}$ & 
& $\left[ \begin{matrix} \cdot & \cdot & x \\ \cdot & x & y \\ x & y & y \end{matrix} \right]_{x,y}$ & \\
\addlinespace[2pt]
$o_{17}$ & $\left[ \begin{matrix} \alpha^{-1}x & y & \cdot \\ y & \beta y - \gamma x & x \\ \cdot & x & y \end{matrix} \right]_{x,y}$ &&& ($**$) \\
\bottomrule
\end{tabular}
\caption{Representatives of line orbits in $\langle \mathcal{V}_3(\F_q) \rangle$ under the action of $K=\PGL(3,\mathbb{F}_q)$ on subspaces of $\PG(\F_q^3 \otimes \F_q^3)$ induced by the action of $\GL(3,\F_q)$ on $3 \times 3$ matrices $M$ given by $M^D = DMD^\top$ (where $D \in \GL(3,\F_q)$). 
Notation is as in Definition~\ref{def:notation}. 
For brevity, the corresponding {\em vector} subspaces $M$ of $\F_q^3 \otimes \F_q^3$ are shown, so that the $K$-orbit representatives themselves are given by $\PG(M)$.
Condition~($*$) is: $v\lambda^2+uv\lambda - 1 \neq 0$ for all $\lambda \in \F_q$, where $v \in \{v_1,v_2\}$ in the case $o_{15}$. 
Condition~($**$) is: $\lambda^3+\gamma \lambda^2- \beta \lambda+ \alpha \neq 0$ for all $\lambda \in \F_q$.
}
\label{table:main}
\end{table}


\subsection*{Tensor orbit $o_4$}
The tensor orbit representative from Table~\ref{pointTable} is $e_1 \otimes e_1 \otimes e_1 + e_2 \otimes e_1 \otimes e_2$. 
Its first contraction space is $M_4 = \langle e_1 \otimes e_1 , e_1 \otimes e_2 \rangle$, and has rank distribution $[q+1,0,0]$. 
Let $B \in \text{GL}(3,\mathbb{F}_q)$ and suppose that the line $\text{PG}(BM_4)$ is contained in $\langle \mathcal{V}_3(\mathbb{F}_q) \rangle$. 
Then, in particular, $\text{PG}(BM_4)$ is contained in the Veronese variety $\cV_3(\F)$, 
a contradiction. 
Therefore, the tensor orbit $o_4$ does not give rise to any line with a symmetric representative.

\subsection*{Tensor orbit $o_5$}
Here we have tensor orbit representative $e_1 \otimes e_1 \otimes e_1 + e_2 \otimes e_2 \otimes e_2$. 
The first contraction space is $M_5 = \langle e_1 \otimes e_1 , e_2 \otimes e_2 \rangle$, with rank distribution $[2,q-1,0]$. 
Note that $\PG(M_5)$ is contained in $\langle \mathcal{V}_3(\mathbb{F}_q) \rangle$: it is the $K$-line orbit representative given in Table~\ref{table:main} (in the second column). 

Given $B\in \GL(3,\mathbb{F}_q)$ and $i\in \{1,2,3\}$, let $B^i\in \F^3_q$ denote the $i$th column vector of $B$.
If $\text{PG}(BM_5)$ is contained in $\langle \mathcal{V}_3(\mathbb{F}_q) \rangle$ for some $B \in \text{GL}(3,\mathbb{F}_q)$, then both of the matrices $B^1 \otimes e_1$ and $B^2 \otimes e_2$ must be symmetric and of rank~$1$. 
This forces $B^1 = \alpha e_1$ and $B^2 = \beta e_2$ for some $\alpha,\beta \in \mathbb{F}_q^\# = \mathbb{F}_q \setminus \{0\}$, and so
\[
B=\left[ \begin{matrix} \alpha & \cdot & * \\ \cdot & \beta & * \\ \cdot & \cdot & * \end{matrix} \right],
\]
where $*$ denotes an unspecified element of $\F_q$ (and $\cdot$ denotes $0$, as per Definition~\ref{def:notation}). 
Therefore, 
\[
BM_5 = \left[ \begin{matrix} \alpha & \cdot & * \\ \cdot & \beta & * \\ \cdot & \cdot & * \end{matrix} \right]
\left[ \begin{matrix} x & \cdot & \cdot \\ \cdot & y & \cdot \\ \cdot & \cdot & \cdot \end{matrix} \right]_{x,y} 
= \left[ \begin{matrix} \alpha x & \cdot & \cdot \\ \cdot & \beta y & \cdot \\ \cdot & \cdot & \cdot \end{matrix} \right]_{x,y}.
\]
Since $\alpha x$ and $\beta y$ range over all values in $\mathbb{F}_q$ as $x$ and $y$ do, we may relabel $\alpha x$ as $x$ and $\beta y$ as $y$ to see that $BM_5 = M_5$. 
That is, $\PG(BM_5)$ is contained in $\langle \mathcal{V}_3(\mathbb{F}_q) \rangle$ if and only if $\PG(BM_5)=\PG(M_5)$, and so the orbit containing $\PG(M_5)$ is the {\em only} $K$-line orbit in $\langle \mathcal{V}_3(\mathbb{F}_q) \rangle$ arising from the tensor orbit $o_5$.

\subsection*{Tensor orbit $o_6$}
This tensor orbit has representative $e_1 \otimes e_1 \otimes e_1 + e_2 \otimes ( e_1 \otimes e_2 + e_2 \otimes e_1 )$. 
The first contraction space is $M_6 = \langle e_1 \otimes e_1 , e_1 \otimes e_2 + e_2 \otimes e_1 \rangle$, with rank distribution $[1,q,0]$. 
As in the previous case, we note that $\PG(M_6)$ is contained in $\langle \mathcal{V}_3(\mathbb{F}_q) \rangle$, and is the representative given in Table~\ref{table:main}. 
Now suppose that $\PG(BM_6)$ is contained in $\langle \mathcal{V}_3(\mathbb{F}_q) \rangle$ for some $B \in \text{GL}(3,\mathbb{F}_q)$. 
Then $B^1 \otimes e_1$ must be symmetric and of rank~$1$, so $B^1 = \alpha e_1$ for some $\alpha \in \mathbb{F}_q^\#$. 
Moreover, $B^1 \otimes e_2 + B^2 \otimes e_1 = \alpha e_1 \otimes e_2 + B^2 \otimes e_1$ must be symmetric and of rank~$2$, so $B^2 = \beta e_1 + \alpha e_2$ for some $\beta \in \mathbb{F}_q$. 
Therefore,
\[
BM_6 = \left[ \begin{matrix} \alpha & \beta & * \\ \cdot & \alpha & * \\ \cdot & \cdot & * \end{matrix} \right]
\left[ \begin{matrix} x & y & \cdot \\ y & \cdot & \cdot \\ \cdot & \cdot & \cdot \end{matrix} \right]_{x,y}  
= \left[ \begin{matrix} \alpha x + \beta y & \alpha y & \cdot \\ \alpha y & \cdot & \cdot \\ \cdot & \cdot & \cdot \end{matrix} \right]_{x,y}.
\]
Again, we may relabel $\alpha x$ as $x$ and $\beta y$ as $y$ to deduce that $BM_6=M_6$. 
Hence, the orbit containing $\PG(M_6)$ is the only $K$-line orbit in $\langle \mathcal{V}_3(\mathbb{F}) \rangle$ arising from the tensor orbit $o_6$.

\subsection*{Tensor orbit $o_7$}
This tensor orbit has representative $e_1 \otimes e_1 \otimes e_3 + e_2 \otimes ( e_1 \otimes e_1 + e_2 \otimes e_2 )$. 
The first contraction space is $M_7 = \langle e_1 \otimes e_3 , e_1 \otimes e_1 + e_2 \otimes e_2 \rangle$, with rank distribution $[1,q,0]$. 
We claim that $\PG(BM_7)$ is not contained in $\langle \mathcal{V}_3(\mathbb{F}_q) \rangle$ for any $B \in \text{GL}(3,\mathbb{F}_q)$. 
If it were, then $B^1 \otimes e_3$ would have to be symmetric and of rank~$1$, forcing $B^1 = \alpha e_3$ for some $\alpha \in \mathbb{F}_q^\#$. 
However, then $B^1 \otimes e_1 + B^2 \otimes e_2 = \alpha e_3 \otimes e_1 + B^2 \otimes e_2$ would not be symmetric, a contradiction.
Hence, the tensor orbit $o_7$ does not give rise to any line with a symmetric representative.

\subsection*{Tensor orbit $o_8$}
Here the tensor orbit representative is $e_1 \otimes e_1 \otimes e_1 + e_2 \otimes ( e_2 \otimes e_2 + e_3 \otimes e_3 )$. 
The first contraction space is $M_8 = \langle e_1 \otimes e_1 , e_2 \otimes e_2 + e_3 \otimes e_3 \rangle$, with rank distribution $[1,1,q-1]$. 
We show that this yields two $K$-orbits of lines in $\langle \mathcal{V}_3(\mathbb{F}_q) \rangle$. 

Suppose that $\PG(BM_8)$ is contained in $\langle \mathcal{V}_3(\mathbb{F}_q) \rangle$ for some $B \in \text{GL}(3,\mathbb{F}_q)$. 
Then $B(e_1 \otimes e_1) = B^1 \otimes e_1$ must be symmetric and of rank~$1$, forcing $B^1 = \alpha e_1$ for some $\alpha \in \mathbb{F}_q^\#$. 
Moreover, $B^2 \otimes e_2 + B^3 \otimes e_3$ must be symmetric and of rank~$2$, so $B^2 = \beta_2 e_2 + \gamma_2 e_3$ and $B^3 = \beta_3 e_2 + \gamma_3 e_3$, with $\gamma_2=\beta_3$ and $\beta_2\gamma_3 - \beta_3^2 = \beta_2 \gamma_3 - \gamma_2 \beta_3 \neq 0$. 
Therefore,
\[
BM_8 = \left[ \begin{matrix} \alpha & \cdot & \cdot \\ \cdot & \beta_2 & \beta_3 \\ \cdot & \beta_3 & \gamma_3 \end{matrix} \right]
\left[ \begin{matrix} x & \cdot & \cdot \\ \cdot & y& \cdot \\ \cdot & \cdot & y \end{matrix} \right]_{x,y} 
= \left[ \begin{matrix} \alpha x &\cdot&\cdot \\ \cdot& \beta_2 y & \beta_3 y \\ \cdot& \beta_3 y & \gamma_3 y \end{matrix} \right]_{x,y}.
\]
We now claim that $\PG(BM_8)$ is $K$-equivalent to either
\begin{equation} \label{o8-1}
\PG \left( \left[ \begin{matrix} \alpha x&\cdot&\cdot \\ \cdot&\cdot&\beta y \\ \cdot&\beta y&\cdot \end{matrix} \right]_{x,y} \right)
\quad \text{or} \quad 
\PG \left( \left[ \begin{matrix} \alpha x&\cdot&\cdot \\ \cdot&\beta y&\cdot \\ \cdot&\cdot&\gamma y \end{matrix} \right]_{x,y} \right) 
\quad \text{for some} \quad \beta,\gamma \in \mathbb{F}_q^\#,
\end{equation}
according to whether $\beta_2=\gamma_3 = 0$ or not. 
This is clear in the case where $\beta_2=\gamma_3 = 0$, as we simply relabel $\beta_3$ as $\beta$. 
On the other hand, if $\gamma_3 \neq 0$ then 
\[
D \left[ \begin{matrix} \alpha x &\cdot&\cdot \\ \cdot& \beta_2 y & \beta_3 y \\ \cdot& \beta_3 y & \gamma_3 y \end{matrix} \right]_{x,y} D^\top = \left[ \begin{matrix} \alpha x&\cdot&\cdot \\ \cdot & (\beta_2-\beta_3^2\gamma_3^{-1})y & \cdot \\ \cdot&\cdot&\gamma_3 y \end{matrix} \right]_{x,y}, 
\quad \text{where} \quad 
D = \left[ \begin{matrix} 1 & \cdot & \cdot \\ \cdot & 1 & -\beta_3\gamma_3^{-1} \\ \cdot & \cdot & 1 \end{matrix} \right];
\]
and if $\beta_2 \neq 0$ then
\[
D \left[ \begin{matrix} \alpha x &\cdot&\cdot \\ \cdot& \beta_2 y & \beta_3 y \\ \cdot& \beta_3 y & \gamma_3 y \end{matrix} \right]_{x,y} D^\top = \left[ \begin{matrix} \alpha x&\cdot&\cdot \\ \cdot & \beta_2 y & \cdot \\ \cdot&\cdot& (\gamma_3-\beta_3^2\beta_2^{-1})y \end{matrix} \right]_{x,y},
\quad \text{where} \quad 
D = \left[ \begin{matrix} 1 & \cdot & \cdot \\ \cdot & 1 & \cdot \\ \cdot & -\beta_3\beta_2^{-1} & 1 \end{matrix} \right]. 
\]
The claim follows upon appropriately relabelling the variables $x$ and $y$. 
By relabelling the constants in \eqref{o8-1}, we then see that $\PG(BM_8)$ is $K$-equivalent to one of the lines
\[
L = \PG\left( \left[ \begin{matrix}  x&\cdot&\cdot \\ \cdot&\cdot& y \\ \cdot& y&\cdot \end{matrix} \right]_{x,y} \right) \quad \text{or} \quad 
L_\gamma = \PG\left( \left[ \begin{matrix}  x&\cdot&\cdot \\ \cdot& y&\cdot \\ \cdot&\cdot&\gamma y \end{matrix} \right]_{x,y} \right) \text{ with } \gamma \in \mathbb{F}_q^\#. 
\]

We now show that $L$ and $L_\gamma$ represent the same $K$-orbit if and only if $q$ is odd and $-\gamma \in \Box$ (that is, $-\gamma$ is a square). 
Disregarding the first row and column, we have
\[
\left[ \begin{matrix} a&b \\ c&d \end{matrix} \right] \left[ \begin{matrix} 0&1 \\ 1&0 \end{matrix} \right] \left[ \begin{matrix} a&c \\ b&d \end{matrix} \right]= 
\left[ \begin{matrix} a&b \\ c&d \end{matrix} \right]\left[ \begin{matrix} b&d \\ a&c \end{matrix} \right] = 
\left[ \begin{matrix} 2ab&ad+bc \\ ad+bc&2cd \end{matrix} \right]=\left[ \begin{matrix} 1&0 \\ 0&\gamma  \end{matrix} \right]
\]
if and only if $2ab=1$, $ad+bc=0$ and $2cd=\gamma$, a contradiction if $q$ is even. 
If $q$ is odd then $b=(2a)^{-1}$, $d=\gamma (2c)^{-1}$ and $-\gamma = (ca^{-1})^2$ for $a\neq 0 \neq c$, so $L$ and $L_\gamma$ represent the same $K$-orbit if and only if $-\gamma \in \Box$. 
Finally, if $q$ is odd and $-1\in \Box$ then $L$ represents the same $K$-orbit as $L_1$, whereas if $-1\notin \Box$ then $L$ represents the same $K$-orbit as $L_\gamma$ with $\gamma\notin \Box$. 
Hence, we have the following two cases, as per Table~\ref{table:main}: if $q$ is even then there are two $K$-orbits, represented by $L$ and $L_1$; 
and if $q$ is odd then the two $K$-orbits are represented by $L_1$ and $L_\gamma$ with $\gamma\notin \Box$.

\subsection*{Tensor orbit $o_9$}
This tensor orbit has representative 
\[
e_1 \otimes e_3 \otimes e_1 + e_2 \otimes ( e_1 \otimes e_1 + e_2 \otimes e_2 + e_3 \otimes e_3 ).
\] 
The first contraction space is $M_9 = \langle e_3 \otimes e_1 , e_1 \otimes e_1 + e_2 \otimes e_2 + e_3 \otimes e_3 \rangle$, with rank distribution $[1,0,q]$. 
If $\PG(BM_9)$ is contained in $\langle \mathcal{V}_3(\mathbb{F}_q) \rangle$ then $B^3 \otimes e_1$ must be symmetric and of rank~$1$, forcing $B^3 = \alpha e_1$ for some $\alpha \in \mathbb{F}_q^\#$. 
Hence, $B^1 \otimes e_1 + B^2 \otimes e_2 + B^3 \otimes e_3 = B^1 \otimes e_1 + B^2 \otimes e_2 + \alpha e_1 \otimes e_3$ must be symmetric and of rank~$3$. 
Writing $B^1 = \alpha_1 e_1 + \beta_1 e_2 + \gamma_1 e_3$ and $B^2 = \alpha_2 e_1 + \beta_2 e_2 + \gamma_2 e_3$, it follows that we must have $\alpha_1=\gamma_2=0$, $\gamma_1=\alpha$, $\alpha_2=\beta_1$ and $\beta_2 \neq 0$. 
Therefore,
\[
BM_9 = \left[ \begin{matrix} \cdot & \beta_1 & \alpha \\ \beta_1 & \beta_2 & \cdot \\ \alpha & \cdot & \cdot \end{matrix} \right]
\left[ \begin{matrix} y & \cdot & \cdot \\ \cdot & y & \cdot \\ x & \cdot & y \end{matrix} \right]_{x,y} 
= \left[ \begin{matrix} \alpha x & \beta_1 y & \alpha y \\ \beta_1 y & \beta_2 y & \cdot \\ \alpha y & \cdot & \cdot \end{matrix} \right]_{x,y}.
\]
Since $\beta_2 \neq 0$, we may set $\beta_2=1$. 
Relabelling also $\beta_1$ as $\beta$ and $x$ as $\alpha^{-1} x$ yields
\[
BM_9 = \left[ \begin{matrix} x & \beta y & \alpha y \\ \beta y & y & \cdot \\ \alpha y & \cdot & \cdot \end{matrix} \right]_{x,y}.
\]
We now see that $\PG(BM_9)$ is $K$-equivalent to the representative shown in Table~\ref{table:main}, because
\[
D \left[ \begin{matrix} x & \beta y & \alpha y \\ \beta y & y & \cdot \\ \alpha y & \cdot & \cdot \end{matrix} \right]_{x,y} D^\top = \left[ \begin{matrix} x & \cdot & y \\ \cdot & y & \cdot \\ y & \cdot & \cdot \end{matrix} \right]_{x,y}, 
\quad \text{where} \quad 
D = \left[ \begin{matrix} 1 & \cdot & \cdot \\ \cdot & 1 & -\beta\alpha^{-1} \\ \cdot & \cdot & \alpha^{-1} \end{matrix} \right].
\]

\subsection*{Tensor orbit $o_{10}$}
This tensor orbit has representative 
\[
e_1 \otimes (e_1 \otimes e_1 + e_2 \otimes e_2 + ue_1 \otimes e_2) + e_2 \otimes (e_1 \otimes e_2 + ve_2 \otimes e_1),
\]
where $v\lambda^2 + uv\lambda - 1 \neq 0$ for all $\lambda \in \mathbb{F}$, namely condition~($*$) in Tables~\ref{pointTable} and~\ref{table:main}. 
The first contraction space is $M_{10} = \langle e_1 \otimes e_1 + e_2 \otimes e_2 + ue_1 \otimes e_2 , e_1 \otimes e_2 + ve_2 \otimes e_1 \rangle$, with rank distribution $[0,q+1,0]$.
If we take
\[
B = \left[ \begin{matrix} 1 & \cdot & \cdot \\ u & v^{-1} & \cdot \\ \cdot & \cdot & * \end{matrix} \right],
\]
then
\[
BM_{10} = \left[ \begin{matrix} 1 & 0 & \cdot \\ u & v^{-1} & \cdot \\ \cdot & \cdot & * \end{matrix} \right] 
\left[ \begin{matrix} x & ux+y & \cdot \\ vy & x & \cdot \\ \cdot & \cdot & \cdot \end{matrix} \right]_{x,y} = 
\left[ \begin{matrix} x & ux+y & \cdot \\ ux+y & u^2x+uy+v^{-1}x & \cdot \\ \cdot & \cdot & \cdot \end{matrix} \right]_{x,y},
\]
so $\PG(BM_{10})$ lies in $\langle \mathcal{V}_3(\mathbb{F}_q) \rangle$. 
By relabelling $ux+y$ as $y$ and $v^{-1}x$ as $x$, we see that $\PG(BM_{10})$ is the $K$-line orbit representative given in Table~\ref{table:main}. 
Now, this line is a constant rank-$2$ line of $2 \times 2$ matrices, so is an external line to a conic. 
Since the group of a conic acts transitively on the set of external lines to the conic, there is only one $K$-line orbit in $\langle \mathcal{V}_3(\mathbb{F}_q) \rangle$ arising from the tensor orbit $o_{10}$.

\subsection*{Tensor orbit $o_{11}$}
This tensor orbit has representative 
\[
e_1 \otimes ( e_1 \otimes e_1 + e_2 \otimes e_2 ) + e_2 \otimes ( e_1 \otimes e_2 + e_2 \otimes e_3 ).
\] 
The first contraction space is $M_{11} = \langle e_1 \otimes e_1 + e_2 \otimes e_2 , e_1 \otimes e_2 + e_2 \otimes e_3 \rangle$, with rank distribution $[0,q+1,0]$.
We claim that $\PG(BM_{11})$ is not contained in $\langle \mathcal{V}_3(\mathbb{F}_q) \rangle$ for any $B \in \text{GL}(3,\mathbb{F}_q)$. 
If it were, then $B^1 \otimes e_1 + B^2 \otimes e_2$ would need to be symmetric, so in particular $B^1$ and $B^2$ would lie in the span of $e_1$ and $e_2$. 
However, then $B^1 \otimes e_2 + B^2 \otimes e_3$ would not be symmetric, a contradiction. 
Hence, the tensor orbit $o_{11}$ does not give rise to any line with a symmetric representative.

\subsection*{Tensor orbit $o_{12}$}
Here we have tensor orbit representative 
\[
e_1 \otimes (e_1 \otimes e_1 + e_2 \otimes e_2) + e_2 \otimes (e_1 \otimes e_3 + e_3 \otimes e_2).
\] 
The first contraction space is $M_{12} = \langle e_1 \otimes e_1 + e_2 \otimes e_2 , e_1 \otimes e_3 + e_3 \otimes e_2 \rangle$, with rank distribution $[0,q+1,0]$. 
If $\PG(BM_{12})$ is contained in $\langle \mathcal{V}_3(\mathbb{F}_q) \rangle$ for some $B \in \text{GL}(3,\mathbb{F}_q)$, then $B^1 \otimes e_1 + B^2 \otimes e_2$ must be symmetric and of rank~$2$, forcing $B^1 = \alpha_1 e_1 + \beta_1 e_2$ and $B^2 = \alpha_2 e_1 + \beta_2 e_2$, with $\alpha_2 = \beta_1$ and $\alpha_1 \beta_2 - \beta_1^2 = \alpha_1 \beta_2 - \alpha_2 \beta_1 \neq 0$.
Writing $B^3 = \alpha_3 e_1 + \beta_3 e_2 + \gamma_3 e_3$, we then have
\[
B^1 \otimes e_3 + B^3 \otimes e_2 = \alpha_1 e_1 \otimes e_3 + \beta_1 e_2 \otimes e_3 + \alpha_3 e_1 \otimes e_2 + \beta_3 e_2 \otimes e_2 + \gamma_3 e_3 \otimes e_2,
\]
which must also be symmetric and of rank~$2$, forcing $\alpha_1 = 0$, $\alpha_3 = 0$ and $\gamma_3 = \beta_1$. 
Hence,
\[
BM_{12} = \left[ \begin{matrix} \cdot & \beta_1 & \cdot \\ \beta_1 & \beta_2 & \beta_3 \\ \cdot & \cdot & \beta_1 \end{matrix} \right]
\left[ \begin{matrix} x & \cdot & y \\ \cdot & x & \cdot \\ \cdot & y & \cdot \end{matrix} \right]_{x,y} 
= \left[ \begin{matrix} \cdot & \beta_1 x & \cdot \\ \beta_1 x & \beta_2 x + \beta_3 y & \beta_1 y \\ \cdot & \beta_1 y & \cdot \end{matrix} \right]_{x,y}.
\]
Since $\beta_1 \neq 0$, we can relabel this as
\[
BM_{12} = \left[ \begin{matrix} \cdot & x & \cdot \\ x & \alpha x + \beta y & y \\ \cdot & y & \cdot \end{matrix} \right]_{x,y} 
\quad \text{for some} \quad \alpha,\beta \in \mathbb{F}_q.
\]

If $q$ is odd then $\PG(BM_{12})$ is $K$-equivalent to the line
\[
L = \PG \left( \left[ \begin{matrix} \cdot & x & \cdot \\ x & \cdot & y \\ \cdot & y & \cdot \end{matrix} \right]_{x,y} \right), 
\]
because 
\[
\left[ \begin{matrix} \cdot & x & \cdot \\ x & \cdot & y \\ \cdot & y & \cdot \end{matrix} \right]_{x,y} = 
D \left[ \begin{matrix} \cdot & x & \cdot \\ x & \alpha x + \beta y & y \\ \cdot & y & \cdot \end{matrix} \right]_{x,y} D^\top, 
\quad \text{where} \quad
D = \left[ \begin{matrix} 1 & \cdot & \cdot \\ -\tfrac{\alpha}{2} & 1 & -\tfrac{\beta}{2} \\ \cdot & \cdot & 1 \end{matrix} \right].
\]
Therefore, there is a single $K$-line orbit, with representative $L$, as per Table~\ref{table:main}. 
Now suppose that $q$ is even. 
In this case we claim that $\PG(BM_{12})$ is $K$-equivalent either to $L$ or to the line
\[
L' = \PG \left( \left[ \begin{matrix} \cdot & x & \cdot \\ x & x+y & y \\ \cdot & y & \cdot \end{matrix} \right]_{x,y} \right), 
\]
according to whether $\alpha=\beta=0$ or not. 
These two lines lie in different $K$-orbits, characterised by the intersection of the line with the nucleus plane (see fact~(F5) in Section~\ref{V3props}): $L$ is contained in the nucleus plane, and $L'$ intersects the nucleus plane in a point. 
It remains to prove the claim. 
If $\alpha = \beta = 0$ then $\PG(BM_{12})=L$. 
If $\alpha \neq 0 \neq \beta$ then $\PG(BM_{12})$ is $K$-equivalent to $\PG(M)$ for
\[
M=\left[ \begin{matrix} \cdot & \alpha x & \cdot \\ \alpha x & \alpha x + \beta y & \beta y \\ \cdot & \beta y & \cdot \end{matrix} \right]_{x,y} = 
D \left[ \begin{matrix} \cdot & x & \cdot \\ x & \alpha x + \beta y & y \\ \cdot & y & \cdot \end{matrix} \right]_{x,y} D^\top, 
\quad \text{where} \quad 
D = \left[ \begin{matrix} \alpha & \cdot & \cdot \\ \cdot & 1 & \cdot \\ \cdot & \cdot & \beta \end{matrix} \right];
\]
and $\PG(M)$ is, in turn, $K$-equivalent to $L$, by relabelling $\alpha x$ as $x$ and $\beta y$ as $y$. 
If $\beta = 0$ and $\alpha \neq 0$ then $\PG(BM_{12})$ is $K$-equivalent to $\PG(M')$ for
\[
M'=\left[ \begin{matrix} \cdot & \alpha (x+y) & \cdot \\ \alpha (x+y) & \alpha x & \alpha y \\ \cdot & \alpha y & \cdot \end{matrix} \right]_{x,y} = 
D \left[ \begin{matrix} \cdot & x & \cdot \\ x & \alpha x & y \\ \cdot & y & \cdot \end{matrix} \right]_{x,y} D^\top, 
\quad \text{where} \quad 
D = \left[ \begin{matrix} \alpha & \cdot & \alpha \\ \cdot & 1 & \cdot \\ \cdot & \cdot & \alpha \end{matrix} \right];
\]
and $\PG(M')$ is also $K$-equivalent to $L$, by relabelling $\alpha(x+y)$ as $x$ and $\alpha y$ as $y$. 
The case where $\alpha=0$ and $\beta \neq 0$ is analogous, and so the proof of the claim is complete.

\subsection*{Tensor orbit $o_{13}$}
This tensor orbit has representative 
\[
e_1 \otimes (e_1 \otimes e_1 + e_2 \otimes e_2) + e_2 \otimes (e_1 \otimes e_2 + e_3 \otimes e_3).
\] 
The first contraction space is $M_{13} = \langle e_1 \otimes e_1 + e_2 \otimes e_2 , e_1 \otimes e_2 + e_3 \otimes e_3 \rangle$, with rank distribution $[0,2,q-1]$. 
Suppose that $\PG(BM_{13})$ is in $\langle \mathcal{V}_3(\mathbb{F}_q) \rangle$ for some $B \in \text{GL}(3,\mathbb{F}_q)$. 
Then $B^1 \otimes e_1 + B^2 \otimes e_2$ must be symmetric and of rank~$2$, so $B^1 = \alpha_1 e_1 + \beta_1 e_2$ and $B^2 = \alpha_2 e_1 + \beta_2 e_2$, with $\alpha_2 = \beta_1$ and $\alpha_1 \beta_2 - \beta_1^2 = \alpha_1 \beta_2 - \alpha_2 \beta_1 \neq 0$.
Writing $B^3 = \alpha_3 e_1 + \beta_3 e_2 + \gamma_3 e_3$, we then have
\[
B^1 \otimes e_2 + B^3 \otimes e_3 = \alpha_1 e_1 \otimes e_2 + \beta_1 e_2 \otimes e_2 + \alpha_3 e_1 \otimes e_3 + \beta_3 e_2 \otimes e_3 + \gamma_3 e_3 \otimes e_3,
\]
which must also be symmetric and of rank~$2$, forcing $\alpha_1 = \alpha_3 = \beta_3 = 0$.
Hence,
\[
BM_{13} = \left[ \begin{matrix} \cdot & \beta_1 & \cdot \\ \beta_1 & \beta_2 & \cdot \\ \cdot & \cdot & \gamma_3 \end{matrix} \right]
\left[ \begin{matrix} x & y & \cdot \\ \cdot & x & \cdot \\ \cdot & \cdot & y \end{matrix} \right]_{x,y} 
= \left[ \begin{matrix} \cdot & \beta_1 x & \cdot \\ \beta_1 x & \beta_2 x + \beta_1 y & \cdot \\ \cdot & \cdot & \gamma_3 y \end{matrix} \right]_{x,y},
\]
or, equivalently,
\[
BM_{13} = \left[ \begin{matrix} \cdot & x & \cdot \\ x & \alpha x + y & \cdot \\ \cdot & \cdot & \gamma y \end{matrix} \right]_{x,y} 
\quad \text{for some} \quad \alpha \in \mathbb{F}_q, \gamma \in \mathbb{F}_q^\#.
\]

First suppose that $q$ is odd. 
Then $\PG(BM_{13})$ is $K$-equivalent to
\[
L_\gamma = \PG \left( \left[ \begin{matrix} \cdot & x & \cdot \\ x & y & \cdot \\ \cdot & \cdot & \gamma y \end{matrix} \right]_{x,y} \right),
\]
because
\[
\left[ \begin{matrix} \cdot & x & \cdot \\ x & y & \cdot \\ \cdot & \cdot & \gamma y \end{matrix} \right]_{x,y} = 
D \left[ \begin{matrix} \cdot & x & \cdot \\ x & \alpha x + y & \cdot \\ \cdot & \cdot & \gamma y \end{matrix} \right]_{x,y} D^\top, 
\quad \text{where} \quad 
D = \left[ \begin{matrix} 1 & \cdot & \cdot \\ -\frac{\alpha}{2} & 1 & \cdot \\ \cdot & \cdot & 1 \end{matrix} \right].
\]
We claim that the lines $L_\gamma$ comprise two $K$-orbits, characterised by whether $\gamma \in \Box$ or not, as indicated in Table~\ref{table:main}.
If $\gamma \in \Box$ then $L_\gamma$ is $K$-equivalent to $L_1$ because $L_1 = D L_\gamma D^\top$ for $D = \text{diag}(1,1,\delta^{-1})$ with $\delta^2=\gamma$. 
If $\gamma, \gamma' \not \in \Box$ then we may write $\gamma' = \gamma \mu^2$ for some $\mu \in \mathbb{F}_q^\#$, and so $L_{\gamma'} = D L_\gamma D^\top$ for $D = \text{diag}(1,1,\mu)$. 
It remains to show that if $\gamma \not \in \Box$ then $L_\gamma$ is not $K$-equivalent to $L_1$. 
To see this, first consider the line $L_\delta$ for an arbitrary $\delta \in \mathbb{F}_q^\#$, and let $P_\delta$ denote the rank-$2$ point on $L_\delta$ obtained by setting $x=0$. 
If $-\delta \in \Box$ then $P_\delta$ is an exterior point (see fact~(F3) of Section~\ref{V3props}), and if $-\delta \not \in \Box$ then $P_\delta$ is an interior point. 
Hence, if $-\delta \in \Box$ and $-\delta' \not \in \Box$ for some $\delta' \in \mathbb{F}_q^\#$, then $L_\delta$ and $L_{\delta'}$ are not $K$-equivalent. 
We now use this observation to show that $L_1$ and $L_\gamma$ are not $K$-equivalent if $\gamma \not \in \Box$, by considering separately the cases where $-1 \in \Box$ and $-1 \not \in \Box$. 
If $-1 \in \Box$ then the point $P_1$ on $L_1$ is an exterior point, but $-\gamma \not \in \Box$ since $\gamma \not \in \Box$ and $-1 \in \Box$, so the point $P_\gamma$ on $L_\gamma$ is an interior point. 
Similarly, if $-1 \not \in \Box$ then $P_1$ is an interior point, but $-\gamma \in \Box$ since $\gamma \not \in \Box$ and $-1 \not \in \Box$, and hence $P_\gamma$ is an exterior point.

Now suppose that $q$ is even. 
If $\alpha \neq 0$ then, because $\gamma \in \mathbb{F}_q^\#$ is a square, say $\gamma = \delta^2$, $\PG(BM_{13})$ is $K$-equivalent to $\PG(M)$ for
\[
M=\left[ \begin{matrix} \cdot & \alpha x & \cdot \\ \alpha x & \alpha x + y & \cdot \\ \cdot & \cdot & y \end{matrix} \right]_{x,y} = 
D \left[ \begin{matrix} \cdot & x & \cdot \\ x & \alpha x + y & \cdot \\ \cdot & \cdot & \gamma y \end{matrix} \right]_{x,y} D^\top, 
\quad \text{where} \quad 
D = \left[ \begin{matrix} \alpha & \cdot & \cdot \\ \cdot & 1 & \cdot \\ \cdot & \cdot & \delta^{-1} \end{matrix} \right].
\]
Relabelling $\alpha x$ as $x$, it follows that $\PG(BM_{13})$ is $K$-equivalent to
\[
L = \PG\left( \left[ \begin{matrix} \cdot & x & \cdot \\ x & x + y & \cdot \\ \cdot & \cdot & y \end{matrix} \right]_{x,y} \right).
\]
If, on the other hand, $\alpha=0$, then $\PG(BM_{13})$ is $K$-equivalent to $L_\gamma$ with $\gamma=1$. 
It remains to show that $L$ and $L_1$ are not $K$-equivalent. 
To see this, observe that $L$ does not intersect the nucleus plane, while $L_1$ intersects the nucleus plane in a point.

\subsection*{Tensor orbit $o_{14}$}
This tensor orbit has representative 
\[
e_1 \otimes (e_1 \otimes e_1 + e_2 \otimes e_2) + e_2 \otimes (e_2 \otimes e_2 + e_3 \otimes e_3).
\] 
The first contraction space is $M_{14} = \langle e_1 \otimes e_1 + e_2 \otimes e_2 , e_2 \otimes e_2 + e_3 \otimes e_3 \rangle$, with rank distribution $[0,3,q-2]$. 
We see that $\PG(M_{14})$ is contained in $\langle \mathcal{V}_3(\mathbb{F}_q) \rangle$ (it is the representative in the second column of Table~\ref{table:main}), and that in order for $\PG(BM_{14})$ to be contained in $\langle \mathcal{V}_3(\mathbb{F}_q) \rangle$, we must have $B = \operatorname{diag}(\alpha,\beta,\gamma)$ for some $\alpha,\beta,\gamma \in \mathbb{F}_q^\#$. 
Hence, 
\[
BM_{14} = \left[ \begin{matrix} \alpha x & \cdot & \cdot \\ \cdot & \beta x + \beta y & \cdot \\ \cdot & \cdot & \gamma y \end{matrix} \right]_{x,y}.
\]
Relabelling, we see that $\PG(BM_{14})$ is equal to the line
\[
L_\gamma = \PG\left( \left[ \begin{matrix} x & \cdot & \cdot \\ \cdot & \gamma (x + y) & \cdot \\ \cdot & \cdot & y \end{matrix} \right]_{x,y} \right) 
\quad \text{for some} \quad \gamma \in \mathbb{F}_q^\#.
\]
The rest of the argument is essentially the same as in the $o_{13}$ case. 
If $q$ is even then every line of the form $L_\gamma$ is $K$-equivalent to $L_1$ because every $\gamma \in \mathbb{F}_q^\#$ is a square; that is, $L_1 = D L_\gamma D^\top$ for $D=\text{diag}(1,\delta^{-1},1)$ with $\delta^2=\gamma$. 
If $q$ is odd then there are two $K$-orbits. 
Indeed, fixing some $\gamma \not \in \Box$, we find that $\PG(BM_{14})$ is $K$-equivalent to either $L_1$ or $L_\gamma$, and that $L_\gamma$ is $K$-equivalent to $L_{\gamma'}$ for every $\gamma' \not \in \Box$. 
Moreover, $L_1$ and $L_\gamma$ are not $K$-equivalent, because for an arbitrary $\delta \in \mathbb{F}_q^\#$, the rank-$2$ points on the line $L_\delta$ that correspond to $x=0$ and $y=0$, respectively, are exterior or interior points according to whether $-\delta \in \Box$ or not.

\subsection*{Tensor orbit $o_{15}$}
In this case the tensor orbit representative is
\[
e_1 \otimes (e_1 \otimes e_1 + e_2 \otimes e_2 + e_3 \otimes e_3 + ue_1 \otimes e_2) + 
e_2 \otimes (e_1 \otimes e_2 + ve_2 \otimes e_1),
\]
where $v\lambda^2 + uv\lambda -1 \neq 0$ for all $\lambda \in \mathbb{F}_q$, namely condition~($*$) in Tables~\ref{pointTable} and~\ref{table:main}. 
The first contraction space is
$
M_{15} = \langle e_1 \otimes e_1 + e_2 \otimes e_2 + e_3 \otimes e_3 + ue_1 \otimes e_2 , e_1 \otimes e_2 + ve_2 \otimes e_1 \rangle
$, 
with rank distribution $[0,1,q]$. 
Observe first that $\PG(B'M_{15})$ is contained in $\langle \mathcal{V}_3(\mathbb{F}_q) \rangle$ for
\[
B' = \left[ \begin{matrix} \cdot & 1 & \cdot \\ 1 & \cdot & \cdot \\ \cdot & \cdot & 1 \end{matrix} \right].
\]
Let us therefore relabel $B'M_{15}$ as $M_{15}$, and also relabel $(u,v)$ as $(s,t)$, so that
\[
M_{15} = \left[ \begin{matrix} ty & x & \cdot \\ x & sx+y & \cdot \\ \cdot & \cdot & x \end{matrix} \right]_{x,y}.
\]
Arguing as in previous cases, we find that in order for $\PG(BM_{15})$ to be contained in $\langle \mathcal{V}_3(\mathbb{F}) \rangle$, we must have
\[
B = \left[ \begin{matrix} \alpha & t\beta & \cdot \\ \beta & \alpha+st\beta & \cdot \\ \cdot & \cdot & \gamma \end{matrix} \right] 
\quad \text{for some } \alpha,\beta,\gamma \in \mathbb{F}_q. 
\]
In particular, $\gamma \neq 0$, and since the action is determined up to a non-zero scalar we may put $\gamma=1$. 
This yields $\PG(BM_{15}) = \PG(M)$, where
\[
M = \left[ \begin{matrix} ty' & x' & \cdot \\ x' & sx'+y' & \cdot \\ \cdot & \cdot &  x \end{matrix} \right]_{x,y}, 
\quad \text{with} \quad x' = \beta ty + \alpha x + st \beta x \quad \text{and} \quad y' = \beta x + \alpha y.
\]

For the sake of presentation, let us now formally state (and prove) the following claim. 

\begin{Claim} \label{o15claim}
$\PG(BM_{15})$ is $K$-equivalent to the line
\[
L(u,v) = \PG\left( \left[ \begin{matrix} vy & x & \cdot \\ x & ux+y & \cdot \\ \cdot & \cdot & x \end{matrix} \right]_{x,y} \right) 
\]
for some $u,v \in \mathbb{F}_q$ satisfying condition~$(*)$ in Tables~\ref{pointTable} and~\ref{table:main}.
\end{Claim}

\begin{Prf}
First suppose that $\alpha \neq 0$, and consider the matrix
\[
D_1 = \left[ \begin{matrix} 1 & \cdot & \cdot \\ -\beta\alpha^{-1} & 1 & \cdot \\ \cdot & \cdot & 1 \end{matrix} \right].
\]
Then, with the vector subspace $M$ defined as above, we have
\[
D_1MD_1^\top = \left[ \begin{matrix} ty' & \delta' \alpha^{-1} x & \cdot \\ \delta' \alpha^{-1} x & \delta'\alpha^{-1}(\alpha^{-1}y'+(s-2\beta\alpha^{-1})x) & \cdot \\ \cdot & \cdot &  x \end{matrix} \right]_{x,y}, 
\quad \text{where} \quad \delta' = \alpha^2 + \alpha\beta st - \beta^2t.
\]
Since $t\lambda^2 + st\lambda - 1 \neq 0$ for all $\lambda \in \mathbb{F}_q$, we may take $\lambda=\beta\alpha^{-1}$ to verify that $\delta' \neq 0$, and hence we may write $\delta^{-1} = \delta'\alpha^{-1}$ for some $\delta \in \mathbb{F}_q^\#$. 
Setting $D_2 = \operatorname{diag}(1,\delta,1)$, we therefore have
\[
D_2(D_1MD_1^\top)D_2^\top = \left[ \begin{matrix} ty' & x & \cdot \\ x & \delta(\alpha^{-1}y'+(s-2\beta\alpha^{-1})x) & \cdot \\ \cdot & \cdot &  x \end{matrix} \right]_{x,y}.
\]
If we now write $u=\delta(s-2\beta\alpha^{-1})$ and $v=t\alpha\delta^{-1}$, and relabel $\delta \alpha^{-1} y'$ as $y$ (noting that this is not the same $y$ as above), then we see that $\PG(BM_{15})=\PG(M)$ is $K$-equivalent to $L(u,v)$. 
We also see that condition~($*$) holds: each matrix with $x \neq 0$ must have rank~$3$, so in particular if we set $x=1$ then the $(3,3)$ minor $vy^2+uvy-1$ must be non-zero for each $y \in \mathbb{F}_q$.

Now suppose that $\alpha =0$ and $s \neq 0$. 
Then
\[
D_1MD_1^\top = \left[ \begin{matrix} s^2t\beta x & st\beta y & \cdot \\ st\beta y & \beta x - st\beta y & \cdot \\ \cdot & \cdot & x \end{matrix} \right]_{x,y}, 
\quad \text{where} \quad
D_1 = \left[ \begin{matrix} s & \cdot & \cdot \\ -s & 1 & \cdot \\ \cdot & \cdot & 1 \end{matrix} \right].
\]
Writing $\delta = s^2t$ and relabelling $st\beta y$ as $y$ then shows that $\PG(BM_{15})$ is $K$-equivalent to $\PG(M')$, where
\[
M' = \left[ \begin{matrix} \beta\delta x & y & \cdot \\ y & \beta x - y & \cdot \\ \cdot & \cdot & x \end{matrix} \right]_{x,y}.
\]
Noting that $\beta \neq 0$ (because $B$ is non-singular) and taking
\[
D_2 = \left[ \begin{matrix} \cdot & 1 & \cdot \\ \beta^{-1} & \beta^{-1} & \cdot \\ \cdot & \cdot & 1 \end{matrix} \right],
\]
we now find that
\[
D_2M'D_2^\top = \left[ \begin{matrix} -\beta^2 y' & x & \cdot \\ x & \beta^{-1}(\delta+2)x+y' & \cdot \\ \cdot & \cdot & x \end{matrix} \right]_{x,y}, \quad \text{where} \quad y'=\beta^2 y' + \beta x, 
\]
and relabelling $y'$ as $y$ shows that $\PG(BM_{15})$ is $K$-equivalent to $L(u,v)$ with $u=\beta^{-1}(\delta+2)$ and $v=-\beta^2$.

Finally, suppose that $s=\alpha=0$. 
Then relabelling $\beta t y$ as $y$ gives
\[
M = \left[ \begin{matrix} t\beta x & y & \cdot \\ y & \beta x & \cdot \\ \cdot & \cdot & x \end{matrix} \right]_{x,y}.
\]
Moreover, $q$ must be odd, because if $s=0$ then $t\lambda^2 \neq 1$ for all $\lambda \in \mathbb{F}_q$, and if $q$ is even then we can take $\lambda$ such that $\lambda^2=t^{-1}$ to yield a contradiction. 
Similarly, $t \neq 1$, because having $t=1$ would imply that $\lambda^2 \neq 1$ for all $\lambda \in \mathbb{F}_q$. 
Now, we have $D_1MD_1^\top = M''$, where
\[
M'' = \left[ \begin{matrix} \beta(1+t) x - 2y & \beta(1-t)x & \cdot \\ \beta(1-t)x & \beta(1+t) x + 2y & \cdot \\ \cdot & \cdot & x \end{matrix} \right]_{x,y}
\quad \text{and} \quad 
D_1 = \left[ \begin{matrix} -1 & 1 & \cdot \\ 1 & 1 & \cdot \\ \cdot & \cdot & 1 \end{matrix} \right].
\]
Since $t \neq 1$ and $\beta \neq 0$, we may write $\delta^{-1} = \beta(1-t)$ for some $\delta \in \mathbb{F}_q$, so that
\[
D_2M''D_2^\top = \left[ \begin{matrix} \delta^2(\beta(1+t) x - 2y) & x & \cdot \\ x & \beta(1+t) x + 2y & \cdot \\ \cdot & \cdot & x \end{matrix} \right]_{x,y} 
\quad \text{for} \quad 
D_2 = \left[ \begin{matrix} \delta & \cdot & \cdot \\ \cdot & 1 & \cdot \\ \cdot & \cdot & 1 \end{matrix} \right].
\]
Writing $y' = 2y-\beta(1+t)x$, we therefore see that $\PG(BM_{15})=\PG(M)$ is $K$-equivalent to 
\[
\PG \left( \left[ \begin{matrix} -\delta^2y' & x & \cdot \\ x & 2\beta(1+t) x + y' & \cdot \\ \cdot & \cdot & x \end{matrix} \right]_{x,y} \right).
\]
Relabelling $y'$ as $y$ shows that this is the line $L(u,v)$ with $v=-\delta^2$ and $u=2\beta(1+t)$.
\end{Prf}

Let us now apply Claim~\ref{o15claim}. 
First suppose that $q$ is odd. 
Then there are at least two $K$-orbits of lines of the form $L(u,v)$, because the unique point of rank~$2$ is either exterior or interior depending on whether $-v \in \Box$ or not. 
We claim that there are {\em exactly} two $K$-orbits, as per Table~\ref{table:main}. 
Consider two such lines $\langle P_2,P_3\rangle$ and $\langle P'_2,P'_3\rangle$, where $P_2$ and $P_2'$ are of rank~$2$, $P_3$ and $P'_3$ are of rank~$3$, and $P_2$ and $P'_2$ are either both exterior points or both interior points with respect to the unique conic plane $\pi$ in which they are contained (see facts~(F1) and~(F2) in Section~\ref{V3props}). 
We may assume that these two lines are $L(u,v)$ and $L(u',v')$, where either both $-v,-v' \in \Box$ or both $-v,-v' \not \in \Box$. 
One then sees that both of the points $P_3$ and $P'_3$ lie on a line through the point $P_1$ corresponding to $\langle e_3\otimes e_3\rangle$, and a point of rank~$2$ in the conic plane $\pi$. 
In fact, each of the planes $\langle P_1,P'_2,P'_3\rangle$ and $\langle P_1,P_2,P_3\rangle$ is a plane on the point $P_1$ intersecting $\pi$ in an external line to the conic. 
Call the corresponding external lines $L$ and $L'$. 
Since the stabiliser of $P_1$ and $\pi$ inside $K$ acts transitively on external lines to the conic in $\pi$ (compare with the case $o_{10}$), we may assume that $L=L'$. 
Since the stabiliser of an external line $M$ inside the group of a conic also acts transitively on both the set of interior points on $L$ and on the set of exterior points on $L$, we may also assume that $P_2=P_2'$, and therefore $v=v'$ (these points correspond to $x=0$). 
Restricting our coordinates to the conic plane $\pi$ corresponding to the top--left $2\times 2$ sub-matrix, and verifying that the point with coordinates $(0,1,u')$ lies on the line $L$ with equation $-X_0-uvX_1+vX_2=0$, we obtain $u=u'$. 
We conclude that $(u,v)=(u',v')$, proving the claim.

If $q$ is even then there is only one $K$-orbit. 
The proof is essentially the same as in the $q$ odd case, except that now the stabiliser of an external line $L$ inside the group of a conic acts transitively on the points of $L$.

\subsection*{Tensor orbit $o_{16}$}
This tensor orbit has representative 
\[
e_1 \otimes (e_1 \otimes e_1 + e_2 \otimes e_2 + e_3 \otimes e_3) + e_2 \otimes (e_1 \otimes e_2 + e_2 \otimes e_3).
\] 
The first contraction space is $M_{16} = \langle e_1 \otimes e_1 + e_2 \otimes e_2 + e_3 \otimes e_3, e_1 \otimes e_2 + e_2 \otimes e_3 \rangle$, with rank distribution $[0,1,q]$. 
The line $\PG(BM_{16})$ is contained in $\langle \mathcal{V}_3(\mathbb{F}_q) \rangle$ if and only if
\[
B = \left[ \begin{matrix} \cdot & \cdot & \alpha \\ \cdot & \alpha & \beta \\ \alpha & \beta & \gamma \end{matrix} \right] 
\]
for some $\alpha,\beta,\gamma \in \F_q$ with $\alpha \neq 0$, and this yields
\[
BM_{16} = \left[ \begin{matrix} \cdot & \cdot & \alpha x \\ \cdot & \alpha x & \beta x + \alpha y \\ \alpha x & \beta x + \alpha y & \gamma x + \beta y \end{matrix} \right]_{x,y}.
\]
After a suitable relabelling, we may write
\[
BM_{16} = \left[ \begin{matrix} \cdot & \cdot & x \\ \cdot & x & y \\ x & y & \alpha x + \beta y \end{matrix} \right]_{x,y}.
\]

If $q$ is odd then $\PG(BM_{16})$ is $K$-equivalent to the representative 
\[
L_1 = \PG\left( \left[ \begin{matrix} \cdot & \cdot & x \\ \cdot & x & y \\ x & y & \cdot \end{matrix} \right]_{x,y} \right)
\]
given in Table~\ref{table:main}, because
\[
\left[ \begin{matrix} \cdot & \cdot & x \\ \cdot & x & y \\ x & y & \cdot \end{matrix} \right]_{x,y} = D \left[ \begin{matrix} \cdot & \cdot & x \\ \cdot & x & y \\ x & y & \alpha x + \beta y \end{matrix} \right]_{x,y} D^\top, 
\quad \text{where} \quad
D = \left[ \begin{matrix} 1 & \cdot & \cdot \\ \tfrac{\beta}{2} & 1 & \cdot \\ -\frac{1}{2}\left(\alpha+\frac{\beta^2}{4}\right) & -\frac{\beta}{2} & 1 \end{matrix} \right].
\]
Now suppose that $q$ is even. 
If $\beta=0$ then, because $\alpha \in \Box$, $\PG(BM_{16})$ is $K$-equivalent to
\[
L_1 = \PG\left( \left[ \begin{matrix} \cdot & \cdot & x \\ \cdot & x & y \\ x & y & \cdot \end{matrix} \right]_{x,y} \right).
\]
We now claim that if $\beta \neq 0$, then $\PG(BM_{16})$ is $K$-equivalent to
\[
L_2 = \PG\left( \left[ \begin{matrix} \cdot & \cdot & x \\ \cdot & x & y \\ x & y & y \end{matrix} \right]_{x,y} \right).
\]
To see this, first observe that 
\[
D \left[ \begin{matrix} \cdot & \cdot & x \\ \cdot & x & y \\ x & y & \alpha x + \beta y \end{matrix} \right]_{x,y} D^\top = \left[ \begin{matrix} \cdot & \cdot & \beta^2 x \\ \cdot & \beta^2 x & \beta y \\ \beta^2 x & \beta y & \alpha x + \beta y \end{matrix} \right]_{x,y},
\quad \text{where} \quad
D = \left[ \begin{matrix} \beta^2 & \cdot & \cdot \\ \cdot & \beta & \cdot \\ \cdot & \cdot & 1 \end{matrix} \right],
\]
and then relabel to obtain
\[
BM_{16} = \left[ \begin{matrix} \cdot & \cdot & x \\ \cdot & x & y \\ x & y & \delta x + y \end{matrix} \right]_{x,y}, 
\quad \text{where} \quad \delta=\alpha\beta^{-2}.
\]
If $\delta=0$ then the projective space obtained from the right-hand side above is equal to $L_2$, so assume now that $\delta \neq 0$. 
Since $BM_{16}$ is spanned by 
\[
P = \left[ \begin{matrix} \cdot & \cdot & 1 \\ \cdot & 1 & \cdot \\ 1 & \cdot & \delta \end{matrix} \right] 
\quad \text{and} \quad Q = \left[ \begin{matrix} \cdot & \cdot & \cdot \\ \cdot & \cdot & 1 \\ \cdot & 1 & 1 \end{matrix} \right],
\]
it is also spanned by $Q$ and $P+\delta Q$. 
Therefore, 
\[
BM_{16} = \left[ \begin{matrix} \cdot & \cdot & x \\ \cdot & x & \delta x + y \\ x & \delta x + y & y \end{matrix} \right]_{x,y},
\]
and so $\PG(BM_{16})$ is $K$-equivalent to $L_2$ because
\[
D \left[ \begin{matrix} \cdot & \cdot & x \\ \cdot & x & \delta x + y \\ x & \delta x + y & y \end{matrix} \right]_{x,y} D^\top = 
\left[ \begin{matrix} \cdot & \cdot & x \\ \cdot & x & y \\ x & y & y \end{matrix} \right]_{x,y}, 
\quad \text{where} \quad D = \left[ \begin{matrix} 1 & \cdot & \cdot \\ -\delta & 1 & \cdot \\ \cdot & \cdot & 1 \end{matrix} \right].
\]
This completes the proof of the claim. 
It remains to show that $L_1$ and $L_2$ are not $K$-equivalent. 
To see this, observe that the point of rank $2$ on $L_1$ lies in the nucleus plane of $\cV(\F_q)$, but the point of rank $2$ on $L_2$ does not.




\subsection*{Tensor orbit $o_{17}$}
The lines corresponding to the tensor orbit $o_{17}$ are constant rank-3 lines, that is, they have rank distribution $[0,0,q+1]$. 
We show that there is a unique $K$-orbit of such lines (see Proposition~\ref{prop3.5}). 
Throughout the proof, we refer to some results proved in Section~\ref{sec:stabs}; the arguments used to prove those results do not in turn depend on any of the arguments given here. 
In particular, we need the following lemma, which counts the number of symmetric representatives of the line orbits arising from $o_{13}$ and $o_{14}$.

\begin{La}\label{lem:nrs} 

The total numbers of symmetric representatives of the line orbits corresponding to the tensor orbits $o_{13}$ and $o_{14}$ are, respectively, $\frac{|K|}{q-1}$ and $\frac{|K|}{6}$.
\end{La}
\begin{Prf}
See the arguments for $o_{13}$ and $o_{14}$ in Section \ref{sec:stabs}, where (in particular) the stabilisers inside $K$ of lines arising from these tensor orbits are determined. 
\end{Prf}

Now, the tensor orbit $o_{17}$ has representative 
\[
e_1 \otimes (e_1 \otimes e_1 + e_2 \otimes e_2 + e_3 \otimes e_3) + e_2 \otimes (e_1 \otimes e_2 + e_2 \otimes e_3 + e_3 \otimes (\alpha e_1 + \beta e_2 + \gamma e_3)),
\]
where $\lambda^3 + \gamma \lambda^2 - \beta \lambda + \alpha \neq 0$ for all $\lambda \in \mathbb{F}_q$, namely condition~($**)$ in Tables~\ref{pointTable} and~\ref{table:main}. 
Note that $\alpha \neq 0$, for otherwise taking $\lambda=0$ would violate condition~($**$).
The first contraction space is
\[
M_{17} = \left[ \begin{matrix} x & y & \cdot \\ \cdot & x & y \\ \alpha y & \beta y & x + \gamma y \end{matrix} \right]_{x,y},
\]
with rank distribution $[0,0,q+1]$ (as noted above).
Setting
\[
B = \left[ \begin{matrix} \alpha & \cdot & \cdot \\ \cdot & -\gamma & 1 \\ \cdot & 1 & \cdot \end{matrix} \right] 
\]
gives
\[
BM_{17} = \left[ \begin{matrix} \alpha^{-1}x & y & \cdot \\ y & \beta y - \gamma x & x \\ \cdot & x & y \end{matrix} \right]_{x,y}, 
\]
so that we can take $\PG(BM_{17})$ as a $K$-orbit representative of lines in $\langle \mathcal{V}_3(\mathbb{F}_q) \rangle$, as in Table~\ref{table:main}. 
We now show that there is only one orbit.

\begin{La}\label{lem:1}
The number of constant rank-3 lines in $\langle \cV_3(\F_q)\rangle$ is $\frac{|K|}{3}$.
\end{La}
\begin{Prf}
Assume that $q>2$, and recall that $|K|=|\PGL(3,\F_q)|=q^3(q^3-1)(q^2-1)$. 
Fix a point $P$ of rank~3. 
Through $P$ there are $q^2+q+1$ lines that contain exactly one point of $\cV_3(\F_q)$. 
According to Table~\ref{pointTable}, the remaining lines through $P$ have rank distributions $[0,i,q+1-i]$ for $i\in \{0,1,2,3\}$. 
Let $n_0,n_1,n_2,n_3$ denote the corresponding numbers of such lines through $P$. 
Then
\begin{equation}\label{eqn:1}
\sum_{i=0}^3 n_i=\frac{q^5-1}{q-1}-(q^2+q+1),
\end{equation}
and counting pairs $(R,L)$ where $R$ is a point of rank 2 and $L=\langle P,R\rangle$ is a line disjoint from $\cV_3(\F_q)$, we obtain
\begin{equation}\label{eqn:2}
\sum_{i=1}^3 in_i=q^2(q^2+q),
\end{equation}
because by Lemma~\ref{lem:counting} there are $q^2$ points of rank 2 on the lines through $P$ and a point of $\cV_3(\F_q)$. 
The lines with rank distribution $[0,2,q-1]$ are symmetric representatives of lines arising from the tensor orbit $o_{13}$. 
By Lemma~\ref{lem:nrs}, there are $\frac{|K|}{q-1}$ such representatives. 
Let $\mathcal{P}_3$ denote the set of rank-3 points in $\langle \cV_3(\F_q) \rangle$. 
Then $|\mathcal{P}_3|=q^5-q^2$, as noted in Section~\ref{V3props}.
Since $K$ acts transitively on $\mathcal{P}_3$, we have
\[
\frac{|{\mathcal{P}}_3| \cdot n_2}{q-1}=\frac{|K|}{q-1},
\]
which implies that $n_2=q^3-q$. 
The lines with rank distribution $[0,3,q-2]$ are symmetric representatives of lines arising from the tensor orbit $o_{14}$. 
By Lemma~\ref{lem:nrs}, there are $\frac{|K|}{6}$ such representatives. 
Since we are assuming that $q \neq 2$, we have
\[
\frac{|{\mathcal{P}}_3|\cdot n_3}{q-2}=\frac{|K|}{6},
\]
which implies that $n_3=\tfrac{1}{6}(q^3-q)(q-2)$. 
Substituting the expressions for $n_2$ and $n_3$ into (\ref{eqn:1}) and (\ref{eqn:2}), we obtain $n_1 = \tfrac{1}{2}(q^4+q^2+2q)$ and hence $n_0 = \tfrac{1}{3}(q^4+q^3-q^2-q)$. 
The total number of constant rank-3 lines is therefore
\[
\frac{|{\mathcal{P}}_3|\cdot n_0}{q+1}=\frac{q^3(q^3-1)(q^2-1)}{3} = \frac{|K|}{3}.
\]
If $q=2$ then the same proof works except that now there are no symmetric representatives of lines corresponding to the tensor orbit $o_{14}$ passing through $P$; that is, $n_3=0$ in this case. 
\end{Prf}

In Proposition~\ref{lem:2}, we determine the stabiliser inside $K$ of a symmetric representative of a line orbit arising from the tensor orbit $o_{17}$. 
It turns out that this stabiliser has order $3$. 
Together with Lemma~\ref{lem:1}, this implies the following result.

\begin{Pro} \label{prop3.5}
There is a unique $K$-orbit of constant rank-3 lines in $\langle \cV_3(\F_q)\rangle$.
\end{Pro}

\begin{Prf}
Immediate from Lemma~\ref{lem:1} and Proposition~\ref{lem:2}. 
\end{Prf}

This concludes the classification of the $K$-orbits of the symmetric representatives of lines arising from the tensor orbits $o_4, \ldots, o_{17}$. 
The results are summarised in Table~\ref{table:main}. 
The number of $K$-line orbits is $15$ for all finite fields $\F_q$. 
Three tensor orbits, namely $o_4$, $o_7$ and $o_{11}$, do not yield a symmetric line orbit representative, and so these are omitted from the table. 
Four orbits split under the action of $K$ for $q$ even but not for $q$ odd, two orbits split for $q$ odd but not for $q$ even, and two orbits split for both $q$ even and $q$ odd.

\section{Line stabilisers for $\F$ a finite field}\label{sec:stabs}

In this section we compute the stabilisers of each of the $K$-orbits of lines in $\langle \cV_3(\F) \rangle$ determined in Section~\ref{sec:orbits}. 
As in that section, we assume here that $\F$ is a finite field $\F_q$. 
The line stabilisers are shown in Table~\ref{table:stabs}. 
The following (common) notation is used: $E_q$ is an elementary abelian group of order $q$, $C_k$ is a cyclic group of order $k$, and $\text{Sym}_n$ is the symmetric group on $n$ letters. 
Moreover, $A \times B$ denotes the direct product of groups $A$ and $B$, while $A:B$ denotes a split extension of $A$ by $B$, with normal subgroup $A$ and subgroup $B$. 
The corresponding total numbers of symmetric line-orbit representatives are readily deduced, and recorded in Table~\ref{table:orbitsLengths}.

For the cases in which a line contains points of rank~2, we also determine for the $q$ odd case how many rank~2 points are exterior (or interior) points, and for the $q$ even case we determine how many rank~2 points lie in the nucleus plane. (Note that in some cases this information has already been obtained as part of the arguments in Section~\ref{sec:orbits}.)

\begin{table}[!t]
\begin{tabular}{cccc}
\toprule
Tensor & \multicolumn{3}{c}{Stabilisers of line orbits in $\langle \mathcal{V}_3(\mathbb{F}_q) \rangle$ under $K=\PGL(3,\F_q)$} \\
\cline{2-4} orbit & Common orbit (all $q$) & \multicolumn{2}{c}{Additional orbit} \\
\cline{3-4} && $q$ odd & $q$ even \\
\hline
\addlinespace[2pt]
$o_5$ 	& $E_q^2 : C_{q-1}^2 : C_2$ & & \\
\addlinespace[2pt]
		\hline
\addlinespace[2pt]	
$o_6$ 	& $E_q^{1+2}:C_{q-1}^2$ & & \\
\addlinespace[2pt]
		\hline
\addlinespace[2pt]	
$o_8$ 	& $C_{q-1} \times \text{O}^\pm(2,\F_q)$, $q \equiv \pm 1(4)$ & $C_{q-1} \times \text{O}^\mp(2,\F_q)$, $q \equiv \pm 1(4)$ & \\
 		& $E_q \times C_{q-1}$, $q$ even & & $C_{q-1} \times \text{SL}(2,\F_q)$ \\
		\hline
\addlinespace[2pt]		
$o_9$	& $E_q^2:C_{q-1}$ & & \\
\addlinespace[2pt]	
		\hline
\addlinespace[2pt]	
$o_{10}$	& $E_q^2:\text{O}^-(2,\F_q)$ & & \\
\addlinespace[2pt]	
		\hline
\addlinespace[2pt]		
$o_{12}$	& $\text{GL}(2,\F_q)$, $q$ odd & & \\
		& $E_q^2:\text{GL}(2,\F_q)$, $q$ even & & $E_q^2:E_q:C_{q-1}$ \\
\addlinespace[2pt]
		\hline
\addlinespace[2pt]	
$o_{13}$	& $C_{q-1} \times C_2$, $q$ odd & $C_{q-1} \times C_2$ & \\
		& $E_q:C_{q-1}$, $q$ even & & $E_q$ \\
\addlinespace[2pt]	
		\hline		
\addlinespace[2pt]		
$o_{14}$ 	& $C_2^2 : \text{Sym}_3$, $q\equiv 1 (4)$ & $C_2^2 : C_2$, $q\equiv 1 (4)$ & \\
		& $C_2^2 : C_2$, $q\equiv 3 (4)$ & $C_2^2 : \text{Sym}_3$, $q\equiv 3 (4)$ & \\
		& $\text{Sym}_3$, $q$ even & & \\
\addlinespace[2pt]
		\hline
\addlinespace[2pt]
$o_{15}$	& $C_2^2$, $q$ odd & $C_2^2$ & \\
		& $C_2$, $q$ even & & \\
\addlinespace[2pt]
		\hline
\addlinespace[2pt]	
$o_{16}$  & $E_q : C_{q-1}$, $q$ odd & & \\
		& $E_q^2 : C_{q-1}$, $q$ even & & $E_q^2$ \\
\addlinespace[2pt]
		\hline
\addlinespace[2pt]		
$o_{17}$ & $C_3$ & & \\
\bottomrule
\end{tabular}
\caption{The stabilisers of line orbits in $\langle \mathcal{V}_3(\mathbb{F}_q) \rangle$ under $K=\PGL(3,\F_q)$.
The layout of the table is consistent with that of Table~\ref{table:main}, that is, the groups in each column are the stabilisers of the orbit representatives shown in the corresponding column of Table~\ref{table:main}. 
For brevity, we write $q \equiv \pm 1(4)$ to mean $q \equiv \pm 1 \pmod 4$.}
\label{table:stabs}
\end{table}

\begin{table}[!t]
\begin{tabular}{ccccc}
\toprule
Tensor & $\sharp$ symmetric line-orbit & $\sharp$ rank-$2$ & $\sharp$ exterior & $\sharp$ points in \\
orbit & representatives & points & points & nucleus plane \\
\midrule
$o_5$ & $\frac{1}{2}q(q+1)(q^2+q+1)$ & $q-1$ & $\tfrac{q-1}{2}$ & $0$ \\
$o_6$ & $(q+1)(q^2+q+1)$ & $q$ & $q$ & $q$ \\
$o_8$ & $q^4(q^2+q+1)$ & $1$ & $0$ or $1$ ($\dagger$) & $0$ or $1$ ($\dagger$) \\
$o_9$ & $q(q^3-1)(q+1)$ & $0$ & -- & -- \\
$o_{10}$ & $\frac{1}{2}q(q^3-1)$ & $q+1$ & $\tfrac{q+1}{2}$ & $0$ \\
$o_{12}$ & $q^2(q^2+q+1)$ & $q+1$ & $q+1$ & $1$ or $q+1$ ($\dagger$) \\
$o_{13}$ & $q^3(q^3-1)(q+1)$ & $2$ & $1$ or $2$ ($\dagger$) & $0$ or $1$ ($\dagger$) \\
$o_{14}$ & $\frac{1}{6}q^3(q^3-1)(q^2-1)$ & $3$ & $1$ or $3$ ($\dagger$) & $0$ \\
$o_{15}$ & $\frac{1}{2}q^3(q^3-1)(q^2-1)$ & $1$ & $0$ or $1$ ($\dagger$) & $0$ \\
$o_{16}$ & $q^2(q^3-1)(q+1)$ & $1$ & $1$ & $0$ or $1$ ($\dagger$) \\
$o_{17}$ & $\frac{1}{3}q^3(q^3-1)(q^2-1)$ & $0$ & -- & -- \\
\bottomrule
\end{tabular}
\caption{The total number of line orbit representatives in $\langle \cV_3(\F_q) \rangle$ corresponding to each tensor orbit in $\F_q^2 \otimes \F_q^3 \otimes \F_q^3$. 
Also shown is the total number of rank-$2$ points on each line (third column), and the total number of these that are exterior points (for $q$ odd, fourth column) or lie in the nucleus plane of $\cV_3(\F_q)$ (for $q$ even, fifth column). 
In some cases, indicated by ($\dagger$), the data in the fourth and/or fifth columns depends on the orbit and, for $q$ odd, (possibly) the parity of $q$ modulo~$4$; in these cases, we refer the reader to the text for full details.}
\label{table:orbitsLengths}
\end{table}

\subsection*{Tensor orbit $o_5$} 
Here there is a single $K$-orbit, represented by the line
\[
L = \PG\left( \left[ \begin{matrix} x & \cdot & \cdot \\ \cdot & y & \cdot \\ \cdot & \cdot & \cdot \end{matrix} \right]_{x,y} \right),
\]
which has rank distribution $[2,q-1,0]$. 
The two points of rank $1$ on $L$ determine a conic in $\cV_3(\F_3)$, and the stabiliser of $L$ must either fix or swap these points. 
In the preimage under the Veronese map $\nu_3$, this corresponds to the setwise stabiliser of two points, inside the stabiliser of a line of $\text{PG}(2,\F_q)$ in $\text{PGL}(3,\F_q)$. 
This group is isomorphic to $E_q^2:C_{q-1}^2:C_2$ (as shown in the second column of Table~\ref{table:stabs}), and so there are 
\[
\frac{|K|}{2q^2(q-1)^2} = \frac{q^3(q^3-1)(q^2-1)}{2q^2(q-1)^2} = \frac{1}{2}q(q+1)(q^2+q+1)
\]
lines in this $K$-line orbit (as shown in Table~\ref{table:orbitsLengths}). 

If $q$ is odd, then the rank-$2$ points on $L$ comprise $\tfrac{q-1}{2}$ exterior points and $\tfrac{q-1}{2}$ interior points; 
if $q$ is even, then all rank-$2$ points on $L$ lie outside the nucleus plane (see Table~\ref{table:orbitsLengths}).

\subsection*{Tensor orbit $o_6$}
Here there is a single $K$-orbit, represented by the line
\[
L = \PG\left( \left[ \begin{matrix} x & y & \cdot \\ y & \cdot & \cdot \\ \cdot & \cdot & \cdot \end{matrix} \right]_{x,y} \right),
\]
which has rank distribution $[1,q,0]$. 
If we consider the unique point $P$ of rank $1$ on $L$, and any other point on $L$, then these two points determine a conic $\cC$. 
The stabiliser of $L$ in $K$ is isomorphic to the stabiliser of the flag $(\nu_3^{-1}(P),\nu_3^{-1}(\cC))$ in $\text{PG}(2,\F_q)$. 
This group is isomorphic to $E_q^{1+2}:C_{q-1}^2$, where the group $E_q^{1+2}$ has centre $Z \cong E_q$ and $E_q^{1+2}/Z \cong E_q^2$.
In particular, there are $(q+1)(q^2+q+1)$ lines in this orbit.

Since $L$ is a tangent of $\cC$, every point on $L$ different from $P$ is an exterior point for $q$ odd. 
If $q$ is even, then every rank-$2$ point on $L$ lies in the nucleus plane.

\subsection*{Tensor orbit $o_8$}
Lines arising from the tensor orbit $o_{8}$ have rank distribution $[1,1,q-1]$. 
For $q$ odd, the unique point of rank $2$ can be an internal point or an external point to the unique conic that it determines.
The stabiliser of an external point (respectively, internal point) inside the group of a non-degenerate conic has size $2(q-1)$ (respectively, $2(q+1)$).
Let $L$ be a symmetric representative of the tensor orbit $o_8$, let $P_1$ be the point of rank $1$ on $L$, and let $P_2$ be the point of rank $2$.
The stabiliser of $P_1$ inside the pointwise stabiliser of the conic plane $\pi=\langle \cC(P_2)\rangle$ in $K$ has size $q-1$, as it corresponds to the group of homologies with common center and axis in $\PG(2,\F_q)$.
The stabiliser of $L$ is therefore isomorphic to one of $C_{q-1} \times \text{O}^\pm(2,\F_q)$, and has size $2(q-1)(q\mp 1)$. Hence, the total number of lines in the two $K$-orbits arising from the tensor orbit $o_8$ is
\[
\frac{|K|}{2(q+1)(q-1)}+\frac{|K|}{2(q-1)^2}=q^4(q^2+q+1).
\]

For $q$ even, the two orbits are characterised by whether or not the unique point of rank $2$ is the nucleus of the unique conic determined by it. 
The point of rank $2$ corresponding to $(x,y)=(0,1)$ in the orbit representative $L$ shown in the fourth column of Table~\ref{table:main} is the nucleus $N$ of the unique conic $\cC$ determined by $N$.
The stabiliser $K_L$ is therefore equal to the stabiliser of $\cC$ and the unique point $R$ of rank $1$ on $L$. 
Note that $R$ is not on $\cC$. 
In the preimage under the Veronese map, this corresponds to the stabiliser of an anti-flag. 
Hence, $K_L$ is isomorphic to $C_{q-1} \times \text{SL}(2,\F_q)$ and has size 
$$q(q^2-1)(q-1).$$
The other orbit is represented by the line $L_1$ shown in the second column of Table~\ref{table:main}. 
The unique point $P$ of rank $2$ on $L_1$ is not the nucleus of the conic $\cC(P)$ that it determines, and so $K_{L_1}$ is equal to the stabiliser of $\cC(P)$, the unique point of rank $1$ on $L$, and the point $Q$ on $\cC(P)$ obtained by intersecting $\cC(P)$ with the unique tangent to $\cC(P)$ through $P$. 
If an element of $K_{L_1}$ fixes a point $Q'$ on $\cC(P)\setminus\{Q\}$, then it fixes the intersection of the line through $Q'$ and $P$ with $\cC(P)$, and therefore fixes $\cC(P)$ pointwise. 
Since the pointwise stabiliser of $\cC(P)$ inside $K_{L_1}$ corresponds to the group of perspectivities with centre not on the axis in the preimage under the Veronese map, it has size $q-1$. 
This implies that $K_{L_1}$ has size $q(q-1)$.
Specifically, $K_{L_1} \cong E_q \times C_{q-1}$.
The total number of lines in these two orbits is therefore
\[
\frac{|K|}{q(q-1)}+\frac{|K|}{q(q^2-1)(q-1)}=q^4(q^2+q+1),
\]
as in the $q$ odd case.

\subsection*{Tensor orbit $o_9$}
Here there is a unique $K$-orbit, represented by the line
\[
L = \PG(\left( \left[ \begin{matrix} x & \cdot & y \\ \cdot & y & \cdot \\ y & \cdot & \cdot \end{matrix} \right]_{x,y} \right),
\]
which has rank distribution $[1,0,q]$. 
Before we determine the stabiliser of $L$ in $K$, let us introduce some terminology and prove two lemmas. 

\begin{La}\label{lem:counting}
A point of rank $3$ in $\langle \cV_3(\F_q)\rangle$ lies on $q^2$ lines with rank distribution $[1,1,q-1]$ and on $q+1$ lines with rank distribution $[1,0,q]$.
\end{La}

\begin{Prf}
It follows from the above treatment of the $o_8$ orbit that there are in total $q^4(q^2+q+1)$ lines with rank distribution $[1,1,q-1]$. 
Since there is just one orbit $\mathcal{P}_3$ of points of rank 3, each such point is on the same number, say $k$, of such lines. 
Counting pairs $(P,L)$ where $P$ is a point of rank $3$ and $L$ is a line with rank distribution $[1,1,q-1]$ containing $P$, we obtain
\[
|{\mathcal{P}}_3|\cdot k = q^4(q^2+q+1)(q-1),
\]
which implies that $k=q+1$. 
The result follows because $\cV_3(\F_q)$ contains $q^2+q+1$ points and a line through a rank-3 point contains at most one point of $\cV_3(\F_q)$.
\end{Prf}

Given a point $P$ of rank 3 in $\langle \cV_3(\F_q)\rangle$, we denote by ${\mathcal{N}}(P)$ the set of points in $\cV_3(\F_q)$ that together with $P$ span a line without rank-2 points. 
It follows from Lemma~\ref{lem:counting} that $|{\mathcal{N}}(P)|=q+1$. 
In the next lemma, we show that ${\mathcal{N}}(P)$ is a normal rational curve (NRC). 
The definition and properties of NRCs may be found in, for instance, \cite[p.~10]{Harris}.

\begin{La}\label{lem:nrc1}
If $P$ is a point of rank $3$ in $\langle \cV_3(\F_q)\rangle$, then the set ${\mathcal{N}}(P)$ is a NRC for $q$ odd (of degree 4 for $q>3$, and of degree 3 for $q=3$), and a conic if $q$ is even.
\end{La}

\begin{Prf}
Since $K$ acts transitively on rank-3 points, we may assume without loss of generality that $P$ corresponds to the matrix with ones on the anti-diagonal and zeroes everywhere else. 
Suppose first that $q$ is odd, and consider the set $S_P$ of images under the Veronese map of the points on the conic $X_0X_2+\frac{1}{2}X_1^2=0$ in $\PG(2,\F_q)$. 
Then $S_P$ is a NRC, since it is the image of a non-degenerate conic (see fact~(F4) of Section~\ref{V3props}). 
A straightforward calculation shows that each of the lines spanned by $P$ and a point of $S_P$ has rank distribution $[1,0,q]$; that is, ${\mathcal{N}}(P)=S_P$ is a NRC. 
For $q$ even, consider the set $S_P$ of images under the Veronese map of the points on the line $X_1=0$. 
Then $S_P$ is a conic, and a straightforward calculation shows that each of the lines spanned by $P$ and a point of $S_P$ has rank distribution $[1,0,q]$; that is, ${\mathcal{N}}(P)=S_P$ is a conic.
\end{Prf}

Now consider the line $L$ (as defined above), and let $R$ be the unique point of rank 1 on $L$. 
Suppose that $q$ is odd. The stabiliser of $R$ inside $K_L$ acts transitively on the points of rank~3 on $L$, since each point of rank~3 determines a NRC through $R$, each of these NRCs is the image of a non-degenerate conic (see fact~(F4) of Section~\ref{V3props}) passing through $\nu_3^{-1}(R)$, and the stabiliser of $\nu_3^{-1}(R)$ inside $\PGL(3,\F_q)$ acts transitively on the non-degenerate conics containing $\nu_3^{-1}(R)$.
Next, consider a point $P$ of rank~3 on $L$. 
Lemma~\ref{lem:nrc1} implies that the stabiliser of $P$ in $K$ equals the stabiliser of ${\mathcal{N}}(P)$ in $K$, which is isomorphic to $\PGL(2,\F_q)$, since ${\mathcal{N}}(P)$ is the image of a non-degenerate conic under the Veronese map. 
The stabiliser of $P$ inside $K_L$ fixes $P$ and $R=L\cap{\mathcal{N}}(P)$, and is therefore isomorphic to the stabiliser in $\PGL(2,\F_q)$ of a point in $\PG(1,\F_q)$, which is $E_q:C_{q-1}$. 
Explicitly, we have $DLD^\top=L$ if and only if $D \in \text{GL}(3,\F_q)$ has the form
\[
D = \left[ \begin{matrix} d_{11} & d_{12} & d_{13} \\ \cdot & d_{22} & -d_{11}^{-1}d_{12}d_{22} \\ \cdot & \cdot & d_{11}^{-1}d_{22}^2 \end{matrix} \right].
\]
The stabiliser of $L$ in $K=\PGL(3,\F_q)$ therefore has order $q^2(q-1)$ 
and is isomorphic to $E_q^2:C_{q-1}$. 
Indeed, the same is true for $q$ even, as can be seen from the explicit form of $D$ given above.
In particular, the $K$-orbit of $L$ has size $q(q^3-1)(q+1)$ both for $q$ even and for $q$ odd.

\subsection*{Tensor orbit $o_{10}$}
\noindent 
A line $L$ in the orbit $o_{10}$ is a line in a conic plane $\pi$ disjoint from the conic $\cC$ consisting of the rank-$1$ points in $\pi$. 
It follows that $\pi$ and $\cC$ are fixed by $K_L$. 
The pointwise stabiliser inside $K$ of the plane $\pi$ corresponds to the pointwise stabiliser of a line in the projectivity group of $\PG(2,\F_q)$, and has size $q^2(q-1)$. 
The stabiliser in $K_\pi$ of an external line $L$ to a conic $\cC$ has size $2(q+1)$. 
We therefore have $K_L \cong E_q^2:\text{O}^-(2,q)$. 
In particular, $|K_L|=2q^2(q-1)(q+1)$, and so the number of symmetric representatives of lines arising from the tensor orbit $o_{10}$ is $\tfrac{1}{2}q(q^3-1)$. 

Moreover, there are $q+1$ tangents to the conic $\cC$, and every tangent meets the external line $L$. 
Since every point is on zero or two tangents, there are $\tfrac{q+1}{2}$ exterior points on $L$ if $q$ is odd (and the other $\tfrac{q+1}{2}$ points on $L$ are interior points). 
If $q$ is even, then all points on $L$ lie outside of the nucleus plane.

\subsection*{Tensor orbit $o_{12}$}
Here, for every $q$, there is a $K$-orbit represented by the line
\[
L = \PG(M) \quad \text{where} \quad
M = \left[ \begin{matrix} \cdot & x & \cdot \\ x & \cdot & y \\ \cdot & y & \cdot \end{matrix} \right]_{x,y},  
\]
which has rank distribution $[0,q+1,0]$. 
If $q$ is odd then we have $DMD^\top=M$ for $D \in \text{GL}(3,\F_q)$ if and only if $D$ has the form 
\[
D = \left[ \begin{matrix} d_{11} & \cdot & d_{13} \\ \cdot & d_{22} & \cdot \\ d_{31} & \cdot & d_{33} \end{matrix} \right].
\] 
Modulo scalars, these matrices comprise a group 
isomorphic to $\text{GL}(2,\F_q)$, 
and so there are $|K|/(q(q^2-1)(q-1)) = q^2(q^2+q+1)$ lines in this $K$-orbit when $q$ is odd. 
If $q$ is even then $DMD^\top=M$ if and only if
\[
D = \left[ \begin{matrix} d_{11} & \cdot & d_{13} \\ d_{21} & d_{22} & d_{23} \\ d_{31} & \cdot & d_{33} \end{matrix} \right].  
\] 
Here there is no restriction on $d_{21}$ or $d_{23}$, so $K_L \cong E_q^2:\text{GL}(2,\F_q)$ 
and the orbit size is $q^2+q+1$. 

When $q$ is even, all points on $L$ lie in the nucleus plane of $\cV_3(\F_q)$. 
For $q$ odd, all points on $L$ are exterior points:

\begin{Pro}
If $q$ is odd and $L$ is a constant rank-$2$ line in $\langle \cV(\F_q)\rangle$ not contained in a conic plane of $\cV(\F_q)$, then every point on $L$ is an exterior point.
\end{Pro}

\begin{Prf} Let $X_0,X_1,X_2$ (respectively $Y_0,\ldots,Y_5$) denote the homogeneous coordinates in $\PG(2,\F_q)$ (respectively $\PG(5,\F_q)$).
There is a unique orbit of such lines, arising from the tensor orbit $o_{12}$. 
Each point $P_{x,y}$ on $L$ is in the unique conic plane $\langle \nu_3(L_{x,y})\rangle$, where $L_{x,y}$ is the line with equation $yX_0-xX_2=0$ in $\PG(2,\F_q)$. 
The image $\nu_3(L_{x,y})$ is the conic with equation $Y_0Y_3-Y_1^2=0$ in the plane $\pi_{x,y}$ with equation $xY_2-yY_0=xY_4-yY_1=x^2Y_5-y^2Y_0=0$. 
The point $P_{x,y}$ is on the tangents $Y_0=0$ and $Y_3=0$ in $\pi_{x,y}$ to $\cC(P_{x,y})$.
\end{Prf}

If $q$ is even then there is a second $K$-orbit, represented by the line 
\[
L_\text{e} = \PG(M_\text{e}) \quad \text{where} \quad 
M_\text{e} = \left[ \begin{matrix} \cdot & x & \cdot \\ x & x+y & y \\ \cdot & y & \cdot \end{matrix} \right]_{x,y}. 
\]
One may check that $DM_\text{e}D^\top=M_\text{e}$ if and only if $D$ has the form
\[
D = \left[ \begin{matrix} d_{11} & \cdot & d_{22}+d_{33} \\ d_{21} & d_{22} & d_{23} \\ d_{11}+d_{22} & \cdot & d_{33} \end{matrix} \right].  
\] 
The stabiliser of $L_\text{e}$ in $K$ is therefore isomorphic to $E_q^2:E_q:C_{q-1}$, which has order $q^3(q-1)$, so $|L_\text{e}^K| = (q^3-1)(q+1)$. 
Hence, when $q$ is even the total number of lines arising from the tensor orbit $o_{12}$ is
\[
(q^2+q+1) + (q^3-1)(q+1) = q^2(q^2+q+1),
\]
which is the same as in the $q$ odd case. 
A point on $L_\text{e}$ lies in the nucleus plane if and only if $x=y$, so $L_\text{e}$ intersects the nucleus plane in one point.

\subsection*{Tensor orbit $o_{13}$}
First consider the case where $q$ is odd. 
There are two orbits, each represented by 
\[
\PG\left( \left[ \begin{matrix} \cdot & x & \cdot \\ x & y & \cdot \\ \cdot & \cdot & \gamma y \end{matrix} \right]_{x,y} \right) 
\quad \text{for some} \quad \gamma \in \F_q^\#.
\] 
The rank-$2$ point corresponding to $y=0$ is always exterior.
If $\gamma=1$, as in the second column of Table~\ref{table:main}, then the rank-$2$ point corresponding to $x=0$ is exterior if $-1 \in \Box$ (that is, if $q \equiv 1 \pmod 4$) and interior otherwise (if $q \equiv 3 \pmod 4$). 
If $\gamma \not \in \Box$, as in the third column of Table~\ref{table:main}, then the situation is reversed: the rank-$2$ point corresponding to $x=0$ is exterior if $q \equiv 3 \pmod 4$ and interior if $q \equiv 1 \pmod 4$. 
Now, let $L$ be the line spanned by the points $P_1$ and $P_2$ of rank 2 corresponding to $y=0$ and $x=0$, respectively, in the above matrix. 
The conics $\cC(P_1)$ and $\cC(P_2)$ (uniquely determined by $P_1$ and $P_2$) intersect in a point $Q$. 
The point $P_1$ is on the tangent line to $\cC(P_1)$ through $Q$. 
The subgroup of the stabiliser of $\cC(P_1)$ fixing $Q$ and $P_1$ also fixes the other point of $\cC(P_1)$ on a tangent to $P_1$, but acts transitively on the remaining points of $\cC(P_1)$. 
It is isomorphic to $C_{q-1}$. 
On the other hand, the point $P_2$ is on a secant through $Q$. 
The subgroup of the stabiliser of $\cC(P_2)$ fixing $Q$ and $P_2$ has order 2 (this is independent of the choice of $\gamma$). 
The stabiliser of $L$ is therefore isomorphic to $C_{q-1} \times C_2$, and so the total number of lines arising from the tensor orbit $o_{13}$ for $q$ odd is
\[
2 \cdot \frac{|K|}{2(q-1)}=q^3(q^3-1)(q+1).
\]

For $q$ even, the first orbit is represented by the above line with $\gamma=1$. 
The rank-2 point $P_1$ corresponding to $y=0$ lies in the nucleus plane, and the other rank-2 point ($P_2$, say, corresponding to $x=0$) does not.
The point $P_1$ is the nucleus of the conic $\cC(P_1)$ and so the stabiliser of $\cC(P_1)$, $Q$ and $P_1$ is isomorphic to $E_q:C_{q-1}$. 
The stabiliser of $\cC(P_2)$, $Q$ and $P_2$ is trivial since $P_2$ is on a secant through $Q$. 
The other orbit is represented by
\[
\PG\left( \left[ \begin{matrix} \cdot & x & \cdot \\ x & x+y & \cdot \\ \cdot & \cdot & y \end{matrix} \right]_{x,y} \right).
\]
Neither rank-2 point lies in the nucleus plane.
The point $P_3$ corresponding to $y=0$ is not the nucleus of $\cC(P_3)$, but it is on the tangent through $Q$. 
The stabiliser of $\cC(P_3)$, $Q$ and $P_3$ is therefore isomorphic to $E_q$. 
The stabiliser of $\cC(P_2)$, $Q$ and $P_2$ is trivial since $P_2$ is on a secant through $Q$. 
We conclude that the total number of lines arising from the tensor orbit $o_{13}$ for $q$ even is 
\[
\frac{|K|}{q(q-1)} + \frac{|K|}{q}=q^3(q^3-1)(q+1),
\]
as in the $q$ odd case.

\subsection*{Tensor orbit $o_{14}$}
Consider the $K$-orbit represented by the line
\[
L_\gamma = \PG(\left( \left[ \begin{matrix} x & \cdot & \cdot \\ \cdot & \gamma( x+y) & \cdot \\ \cdot & \cdot & y \end{matrix} \right]_{x,y} \right) \quad \text{for some} \quad \gamma \in \F_q^\#,
\]
with rank distribution $[0,3,q-2]$. 

First suppose that $q$ is odd. 
The rank-2 point $P_\text{e}$ obtained for $(x,y)=(1,-1)$ is always an exterior point, while the other two rank-2 points, namely $P_1$ obtained for $(x,y)=(1,0)$ and $P_2$ obtained for $(x,y)=(0,1)$, are both exterior if $-\gamma \in \Box$ and both interior otherwise. 
In particular, if $\gamma=1$, as in the second column of Table~\ref{table:main}, then there are three exterior points if $q \equiv 1 \pmod 4$, and one exterior point if $q \equiv 3 \pmod 4$. 
If $\gamma \not \in \Box$, as in the third column of Table~\ref{table:main}, then the situation is reversed.
Now, the conics $\cC(P_\text{e})$ and $\cC(P_i)$ meet in a point $Q_i$ ($i=1,2$), and the conics $\cC(P_1)$ and $\cC(P_2)$ meet in a point $Q_{12}$. 
For each $i\in\{1,2\}$, the point $P_i$ is on the secant through $Q_i$ and $Q_{12}$, and the point $P_\text{e}$ is on the secant through $Q_1$ and $Q_2$. 
The subgroup of the group of a conic stabilising two points on the conic and a third point on the secant through these two points has order $2$. 
By considering two of the three conic planes, this gives us a group of order $4$. 
By considering the conics in the preimage of the Veronese map, we obtain a triangle, from which one observes that the action on two of the sides determines the action on the third side. 
This implies that the action on two of the conic planes determines the action on the third conic plane. 
We conclude that the subgroup of $K_{L_\gamma}$ stabilising the points $P_\text{e}$, $P_1$ and $P_2$ has order $4$. 
Taking into account the permutations of the points $P_\text{e}$, $P_1$ and $P_2$ in the case where all three points are exterior, this amounts to a group of order $24$, isomorphic to $(C_2 \times C_2) : \text{Sym}_3$. 
In the other case, $K_{L_\gamma}$ is isomorphic to $(C_2 \times C_2) : C_2$ and has order $8$.

Now consider the case where $q$ is even. 
Note that all points of rank 2 on $L_\gamma$ lie outside the nucleus plane.
If $P_1$, $P_2$ and $P_3$ denote the points of rank 2 on $L_\gamma$, then each point $P_i$ is on the secant to the conic $\cC(P_i)$ passing through the intersection points of $\cC(P_i)$ with the other two conics $\cC(P_j)$ and $\cC(P_k)$, where $\{i,j,k\}=\{1,2,3\}$. 
The group fixing the conic $\cC(P_i)$, two points on $\cC(P_i)$ and a point $P_i$ on the secant passing though these two points is trivial, because for $q$ even this group also fixes the unique tangent through $P_i$. 
Taking into account the permutations of the points $P_1$, $P_2$ and $P_3$, we obtain $K_L\cong \text{Sym}_3$ and $|K_L|=6$.

Since $\tfrac{1}{24}+\tfrac{1}{8}=\tfrac{1}{6}$, we conclude that the tensor orbit $o_{14}$ yields $\tfrac{|K|}{6}$ 
symmetric representatives of lines for every $q$.

\subsection*{Tensor orbit $o_{15}$}
Here every $K$-orbit is represented by a line
\[
L = \PG\left( \left[ \begin{matrix} vy & x & \cdot \\ x & ux+y & \cdot \\ \cdot & \cdot & x \end{matrix} \right]_{x,y} \right),
\]
for some $v$. 
The rank distribution is $[0,1,q]$. 
For $q$ odd, the unique point of rank 2 is exterior when $-v \in \Box$ (as in the second column of Table~\ref{table:main}) and interior when $-v \not \in \Box$ (third column); for $q$ even, the unique point of rank 2 lies outside the nucleus plane.
Let $\pi$ denote the plane containing the conic $\cC(R)$ uniquely determined by the point $R$ of rank 2 on $L$.
The group $K_L$ fixing $L$ also fixes the point $P$ in $\cV_3(\F_q)$ corresponding to $e_3\otimes e_3$, and therefore also fixes the line $\ell$ obtained by projecting $L$ from $P$ on to $\pi$. 
This projection $\ell$ corresponds to the $2\times 2$ sub-matrix obtained by deleting the last row and the last column from the above matrix representation of $L$, and is therefore a line through the point $R$ that is external to the conic $\cC(R)$. 
The group $K_L$ is the stabiliser of $P$, $\cC(R)$, $\ell$ and $R$.
The linewise stabiliser of $\ell$ in $K_L$ must fix the set $\{P,P^q\}$ of two conjugate points over the quadratic extension of $\F_q$. 
This implies that the stabiliser $K_L$ must fix $R$ and $\{P,P^q\}$, and must therefore have order twice the order of the pointwise stabiliser of $\ell$ in $K_L$, which has order $2$ for $q$ odd and is trivial for $q$ even. 
We conclude that $K_L \cong C_2^2$ if $q$ is odd, and $K_L \cong C_2$ if $q$ is even.

Since the $G$-line orbit arising from the tensor orbit $o_{15}$ splits into two $K$-orbits for $q$ odd, there are in total $\frac{|K|}{2}$ symmetric representatives of lines for both $q$ even and $q$ odd.

\subsection*{Tensor orbit $o_{16}$}
Here for every $q$, there is a $K$-orbit represented by the line
\[
L = \PG(M) \quad \text{where} \quad M= \left[ \begin{matrix} \cdot & \cdot & x \\ \cdot & x & y \\ x & y & \cdot \end{matrix} \right]_{x,y}. 
\]
The rank distribution is $[0,1,q]$, so if $DMD^\top=M$ for $D = (d_{ij}) \in \text{GL}(3,\F_q)$ then the rank-$2$ point corresponding to $(x,y)=(0,1)$ must be fixed, so $d_{12}=d_{13}=0$ and 
\[
\left[ \begin{matrix} d_{22} & d_{23} \\ d_{32} & d_{33} \end{matrix} \right] 
\left[ \begin{matrix} \cdot & 1 \\ 1 & \cdot \end{matrix} \right]
\left[ \begin{matrix} d_{22} & d_{23} \\ d_{32} & d_{33} \end{matrix} \right]^\top = 
\left[ \begin{matrix} \cdot & \alpha \\ \alpha & \cdot \end{matrix} \right] \quad \text{for some} \quad \alpha \in \F_q^\#.
\]
This makes the $(1,2)$ entry of $DMD^\top$ equal to $d_{11}d_{23}$, which forces $d_{23}=0$ because $D$ must be invertible. 
In particular, $D$ must be lower triangular. 
If $q$ is odd then by considering again the image of the unique rank-$2$ point on $L$, we deduce that $d_{32}=0$. 
By considering the image of an arbitrary point on $L$, we then see that $d_{31}=0$ and $d_{22}^2=d_{11}d_{33}$, so that
\[
D = \left[ \begin{matrix} d_{11} & \cdot & \cdot \\ d_{21} & d_{22} & \cdot \\ \cdot & \cdot & -d_{11}^{-1}d_{22}^2 \end{matrix} \right].
\]
These matrices comprise a group of order $q(q-1)^2$, and upon quotienting out the centre of $\text{GL}(3,\F_q)$ we see that the stabiliser of $L$ in $K$ has order $q(q-1)$, so that the orbit of $L$ has size $q^2(q^3-1)(q+1)$. 
Now suppose that $q$ is even, and consider the image of the point corresponding to $(x,y)=(1,0)$ under a lower-triangular matrix $D$. 
The $(3,3)$ entry is $2d_{31}d_{33}+d_{32}^2=d_{32}^2$, so we again deduce that $d_{32}=0$ (as in the $q$ odd case), but we do {\em not} need $d_{31}=0$. 
We also have $d_{22}^2=d_{11}d_{33}$ in the $q$ even case, so that $D$ has the same form as above, except with no restriction on $d_{31}$.
The stabiliser of $L$ in $K$ therefore has order $q^2(q-1)$, and so the orbit has size $q(q^3-1)(q+1)$. 

If $q$ is even then we also have a second $K$-orbit, represented by the line
\[
L_\text{e} = \PG(M_\text{e}) \quad \text{where} \quad M_\text{e} = \left[ \begin{matrix} \cdot & \cdot & x \\ \cdot & x & y \\ x & y & y \end{matrix} \right]_{x,y}.
\]
If $DM_\text{e}D^\top=M_\text{e}$ for $D = (d_{ij}) \in \text{GL}(3,\F_q)$, then again the rank-$2$ point corresponding to $(x,y)=(0,1)$ must be fixed. 
This forces $D$ to be lower triangular with $d_{33}=d_{22}$. 
By then considering an arbitrary point on $L_\text{e}$, we deduce that $D$ must have the form
\[
D = \left[ \begin{matrix} d_{11} & \cdot & \cdot \\ d_{11}^{-1}d_{32}^2-d_{32} & d_{11} & \cdot \\ d_{31} & d_{32} & d_{11} \end{matrix} \right].
\]
These matrices comprise a subgroup of order $q^2(q-1)$ in $\text{GL}(3,\F_q)$, so the stabiliser of $L_\text{e}$ in $K$ has order $q^2$ and hence the orbit has size $q(q^3-1)(q^2-1)$. 

Therefore, in total there are $q^2(q^3-1)(q+1)$ lines in $\langle \cV(\F_q) \rangle$ arising from the tensor orbit $o_{16}$, whether $q$ is even or odd.
When $q$ is odd, the unique rank-$2$ point is always exterior; when $q$ is even, the unique rank-2 point lies in the nucleus plane for the line $L$ but not for the line $L_\text{e}$.

\subsection*{Tensor orbit $o_{17}$} 
For this final case, we show that the line stabiliser has order $3$. 
Recall that, by Lemma \ref{lem:nrc1}, each point $P$ of rank $3$ defines a NRC ${\mathcal{N}}(P)$ contained in $\cV_3(\F_q)$. 
The following lemma is proved via a straightforward calculation.

\begin{La}\label{lem:polarity}
If $q$ is odd then the map $\rho~:~\cP_3\rightarrow \PG(5,\F_q)$ given by $\rho(P) = \langle {\mathcal{N}}(P)\rangle$ defines the polarity $(y_0,y_1,\ldots,y_5)\mapsto y_0Y_0+y_1Y_1+y_2Y_2+2y_3Y_3+2y_4Y_4+2y_5Y_5=0$ in $\PG(5,\F_q)$.
\end{La}


\begin{La}\label{lem:one_point}
If $P$ and $P'$  are two distinct points on a constant rank-3 line in $\langle \cV_3(\F)\rangle$, then ${\mathcal{N}}(P)$ and ${\mathcal{N}}(P')$ intersect in at most one point. 
\end{La}

\begin{Prf}
If $q$ is even then the statement follows immediately from the fact that each two conics on the quadric Veronesean intersect in a point. 
Now let $q$ be odd and suppose that $W={\mathcal{N}}(P)\cap {\mathcal{N}}(P')$ contains two distinct points $R$ and $Q$. 
Then $\langle \cC(R,Q)\rangle$ intersects $\langle W\rangle$ in at least one line, and there exists a hyperplane through $W$ that contains two conics of the Veronesean. 
Since the map $\rho$ defined in Lemma~\ref{lem:polarity} is a polarity, this hyperplane is the image of a point $S$ on the line through $P$ and $P'$. However, since the hyperplane $S^\rho$ contains two conics, it is not a NRC, and therefore $S$ does not have rank 3, a contradiction.
\end{Prf}

\begin{Pro}\label{lem:2}
The linewise stabiliser in $K$ of a constant rank-$3$ line has order $3$.
\end{Pro}

\begin{Prf}
The linewise stabiliser $K_L$ inside $K$ of a constant rank-3 line $L$ in $\langle \cV_3(\F_q)\rangle$ must fix the set $\{P,P^q,P^{q^2}\}$ of three conjugate points of rank 2 on the line $\overline{L}$ defined over the cubic extension of $\F_q$. 
Also, no element of $K$ can fix one of these three points unless it acts as the identity on the line $\overline{L}$. 
For instance, if $g\in K_L$ fixes $P$, then $g$ must fix $P^q+P^{q^2}$ and $P+P^q+P^{q^2}$, which implies that $g$ fixes a frame of $\overline{L}$ and must therefore fix every point of $\overline{L}$. 
Next we prove that the pointwise stabiliser of $L$ is trivial. 
If $q$ is odd then, by Lemmas~\ref{lem:nrc1} and~\ref{lem:one_point}, any projectivity $\varphi$ fixing $L$ pointwise must fix $q+1$ NRCs pairwise intersecting in a point. 
If $q$ is even then the same lemmas imply that $\varphi$ must fix $q+1$ conics pairwise intersecting in a point. 
In both cases, the set of intersection points contains the image of a frame of $\PG(2,\F_q)$ under the Veronese map, and so $\varphi$ is the identity.
It follows that $K_L$ has order 3.
 \end{Prf}

\section{Algebraically closed fields and the real numbers}\label{sec:otherF}

In this section, we explain how the arguments from the case where $\F$ is a finite field can be modified to treat algebraically closed fields and the case $\F=\mathbb{R}$.
When $\F$ is algebraically closed, the orbits $o_{10}$, $o_{15}$ and $o_{17}$ do not occur in the classification of tensors in $\F^2\otimes \F^3\otimes \F^3$ given in \cite{LaSh2015}, and so in particular we do not obtain the corresponding $K$-line orbits in $\langle \cV_3(\F) \rangle$. 
On the other hand, unlike in that classification, in the study of the symmetric representation of the corresponding line orbits we need to distinguish between the cases $\operatorname{char}(\F)=2$ and $\operatorname{char}(\F)\neq 2$.

\subsection*{Algebraically closed fields $\F$ with $\operatorname{char}(\F)\neq 2$}

The orbits listed in the third column of Table \ref{table:main} (the `additional orbit, $q$ odd' column) do not arise, because these depend on the existence of a non-square in $\F$. 
The orbits in the fourth column also do not arise, because their representatives are $K$-equivalent to the corresponding representatives in the second column for $\operatorname{char}(\F) \neq 2$. 
Hence, the only tensor orbits that yield lines with symmetric representatives are $o_5$, $o_6$, $o_8$, $o_9$, $o_{12}$, $o_{13}$, $o_{14}$ and $o_{16}$, and none of these eight orbits splits under $K$.

\subsection*{Algebraically closed fields $\F$ with $\operatorname{char}(\F)=2$}

In this case the representatives in the fourth column of Table \ref{table:main} {\em do} arise, because they essentially depend on the existence of a nucleus of a non-degenerate conic in $\PG(2,\F)$, a property that holds whenever $\operatorname{char}(\F)=2$. 
We therefore obtain the same eight $K$-orbits from the $\operatorname{char}(\F)\neq 2$ case, plus the four extra $K$-orbits corresponding to the representatives in the fourth column (for $o_8$, $o_{12}$, $o_{13}$ and $o_{16}$).

\subsection*{The real numbers}
Finally, consider the case $\F=\mathbb R$. 
Observe first that the orbit $o_{17}$ does not yield any lines with symmetric representatives, because every cubic polynomial with real coefficients has at least one real root, and so condition~($**$) in Table~\ref{table:main} cannot hold. 
The line orbits corresponding to $o_{8}$, $o_{13}$ and $o_{14}$ split, with representatives as in the second and third columns of Table \ref{table:main}, as the existence of the representatives in the third column depends only on the existence of a non-square $\gamma\in {\mathbb{R}}$ (so one can take $\gamma<0$). 
However, the line orbit corresponding to $o_{15}$ does {\em not} split, because condition~($*$) is equivalent to $u^2v^2+4v=v(u^2v+4)$ being negative, and this implies that $v$ is negative, so the case $-v \not \in \Box$ does not occur. 
In summary, we have a total of 13 $K$-line orbits: one arising from each of the tensor orbits $o_6$, $o_9$, $o_{10}$, $o_{12}$, $o_{15}$, and $o_{16}$, with representatives as in the second column of Table~\ref{table:main}; and two arising from each of $o_8$, $o_{13}$ and $o_{14}$, with representatives as in the second and third columns of Table~\ref{table:main}.

\section{The classification of pencils of conics in $\PG(2,\F_q)$}\label{sec:pencils}

As mentioned in Section~\ref{sec:history}, our results imply the classification of pencils of conics in $\PG(2,\F)$ when $\operatorname{char}(\F) \neq 2$. 
This follows from the following observation.
A pencil of conics is a one-dimensional linear system of quadrics in $\PG(2,\F)$.
If $b$ is the bilinear form associated to the quadratic form $f$ defining a conic $\cC$ in $\PG(2,\F)$, and $B$ is the matrix of $b$ with respect to some basis of $\F^3$, then the conic $\cC$ consists of points whose coordinate vectors $v\in \F^3$ satisfy $vBv^T=0$. 
A projectivity of $\PGL(3,\F)$ induced by the matrix $A\in \GL(3,\F)$ mapping a point with coordinate vector $v$ to the point with coordinate vector $vA$ maps the conic with equation $vBv^T$ to the conic with equation $vABA^Tv^T$. 
Hence, the equivalence classes of pencils of conics under the projectivity group $K=\PGL(3,\F)$ are equivalent to the $K$-orbits of lines in the projective space of symmetric $3\times 3$ matrices.
If $\F$ is algebraically closed with $\operatorname{char}(\F) \neq 2$, it therefore follows from Section~\ref{sec:otherF} that there are 8 equivalence classes of pencils of conics; similarly, there are 13 equivalence classes when $\F=\mathbb{R}$. 
This, of course, agrees with the results of Jordan~\cite{Jordan1906,Jordan1907}. 
For $\F=\F_q$ with $q$ odd, it follows from our classification of $K$-orbits on lines in $\PG(5,\F_q)$ that there are 15 equivalence classes of pencils of conics. 
Representatives of each equivalence class of pencils are given in Table~\ref{pencilTable}. 
This is in agreement with the results of Dickson~\cite{Dickson1908}. 
In addition to Dickson's classification, we have also determined the stabiliser for each equivalence classe, as well as the number of pencils in each class (see Tables~\ref{table:stabs} and~\ref{table:orbitsLengths}). 

\begin{table}[!t]
\begin{tabular}{ccc}
\toprule
Tensor & \multicolumn{2}{c}{Equivalence classes of pencils of conics in $\PG(2,\F_q)$, $q$ odd} \\
\midrule
$o_5$ 	& $(X^2,Y^2)$ & \\		
$o_6$ 	& $(X^2,2XY)$ &  \\		
$o_8$ 	& $(X^2,Y^2+Z^2)$& $(X^2,Y^2+\gamma Z^2)$  \\	
$o_9$	& $(X^2,Y^2+2XZ)$ &  \\		
$o_{10}$	& $(vX^2+Y^2,2XY+uY^2)$ & \\		
$o_{12}$	& $(2XY,2YZ)$ &  \\		
$o_{13}$	& $(2XY,Y^2+Z^2)$ & $(2XY,Y^2+\gamma Z^2)$  \\	
$o_{14}$ 	& $(X^2+Y^2,Y^2+Z^2)$ & $(X^2+\gamma Y^2,\gamma Y^2+Z^2)$  \\	
$o_{15}$	& $(2XY+uY^2+Z^2,v_1X^2+Y^2)$ & $(2XY+uY^2+Z^2,v_2X^2+Y^2)$  \\	
$o_{16}$  & $(2XZ+Y^2,2YZ)$ & \\
$o_{17}$ & $(\frac{1}{\alpha}X^2-\gamma Y^2+2YZ,2XY+\beta Y^2+Z^2)$ &  \\
\bottomrule
\end{tabular}
\caption{
The equivalence classes of pencils of conics in $\PG(2,\F_q)$, $q$ odd, where the parameters $\alpha, \beta, \gamma, u,v,v_1,v_2 \in \F_q$ correspond to those in Table~\ref{table:main}.
}
\label{pencilTable}
\end{table}



\section*{Acknowledgements} 
The first author acknowledges the support of {\em The Scientific and Technological Research Council of Turkey} T\"UB\.{I}TAK (project no.~118F159). 
The second author acknowledges the support of the Australian Research Council Discovery Grant DP140100416, which funded his previous appointment at The University of Western Australia (UWA). 
He is also indebted to UWA's Centre for the Mathematics of Symmetry and Computation for partially funding his visit to the University of Padua in June 2016, during which this work was initiated, and to the University of Padua for their hospitality. 
Both authors would like to thank the anonymous referees for their time and interest in this work, and for pointing out the papers \cite{ArNu2002} and \cite{Wall1977}. 
They would also like to thank Hans Havlicek for valuable discussions and feedback which lead to significant corrections and improvements to the paper.

\end{document}